\documentclass{amsart}
\usepackage{amssymb}
\usepackage{amscd}
\usepackage{eepic}
\usepackage{epic}
\usepackage{hyperref}

\swapnumbers

\newtheorem{Lemma}{Lemma}
\newtheorem{Theorem}[Lemma]{Theorem}
\newtheorem{Proposition}[Lemma]{Proposition}
\newtheorem{Corollary}[Lemma]{Corollary}
\theoremstyle{remark} 
\newtheorem{Remark}[Lemma]{Remark}
\theoremstyle{definition}
\newtheorem{Example}[Lemma]{Example}
\newtheorem{noTitle}[Lemma]{}

\numberwithin{Lemma}{section}

\newcommand{\om}{\omega}
\def \C {{\mathbb C}}
\def \R {{\mathbb R}}
\def \T {{\mathbb T}}
\def \Z {{\mathbb Z}}

\def \half {\frac{1}{2}}
\def \B {{\mathcal{B}}}
\def \J {{\mathcal{J}}}
\def \M {{\mathcal{M}}}

\def \CP {{\mathbb C}{\mathbb P}}
\def \tCP {{\mathbb C}{\mathbb P}{^2}}
\def \calD {{\mathcal D}}
\DeclareMathOperator \Sympl {Sympl}
\DeclareMathOperator \Diff {Diff}

\DeclareMathOperator \GL {GL}
\DeclareMathOperator \PSL {PSL}
\DeclareMathOperator \id {id}
\DeclareMathOperator \Lie {Lie}
\DeclareMathOperator \AGL {AGL}
\DeclareMathOperator \image {image}
\DeclareMathOperator \length {length}
\DeclareMathOperator \perimeter {perimeter}
\DeclareMathOperator \area {area}
\DeclareMathOperator \coker {coker}

\def \calA {{\mathcal A}}
\def \calN {{\mathcal N}}

\def \tM {{\tilde{M}}}

\def \tDelta {{\tilde{\Delta}}}

\def \tgamma {{\tilde{\gamma}}}

\def \Inv {^{-1}}
\def \ssminus {\smallsetminus}
\def \eps {\epsilon}
\def \sss {\scriptstyle}
\def \Cinf {C^\infty}

\def\eoe{\unskip\ \hglue0mm\hfill$\between$\smallskip\goodbreak}

\newcommand     {\comment}[1]   {}
\newcommand{\mute}[2] {}
\newcommand     {\printname}[1] {}
\newcommand{\printversion}{}
\newcommand{\labell}[1] {\label{#1}\printname{#1}}

\begin{document}

\printversion

\title[A compact symplectic four-manifold admits only
finitely many toric actions]{A compact symplectic four-manifold admits only
finitely many inequivalent toric actions}

\author{Yael Karshon}
\address{Department of Mathematics, University of Toronto, 
Toronto, Ontario, Canada, M5S 2E4.}
\email{karshon@math.toronto.edu}

\author{Liat Kessler}
\address{Courant Institute of Mathematical Sciences,
New York University, New York, NY 10012, U.S.A.}
\email{kessler@cims.nyu.edu}

\author{Martin Pinsonnault}
\address{Department of Mathematics, University of Toronto, Toronto, 
Ontario, Canada, M5S 2E4.}
\curraddr{Fields Institute, Toronto, ON, Canada, M5T 3J1.}
\email{mpinsonn@fields.utoronto.ca}

\begin{abstract}
Let $(M,\omega)$ be a four dimensional compact connected symplectic 
manifold.  We prove that $(M,\omega)$ admits only finitely many 
inequivalent Hamiltonian effective 2-torus actions. Consequently, 
if $M$ is simply connected, the number of conjugacy classes 
of 2-tori in the symplectomorphism group $\Sympl(M,\omega)$ is finite. 
Our proof is ``soft''. The proof uses the fact that for symplectic
blow-ups of $\CP^2$ the restriction of the period map to the set of
exceptional homology classes is proper.
In an appendix, we describe results of McDuff that give a properness
result for a general compact symplectic four-manifold, using the theory
of J-holomorphic curves. 
\end{abstract}

\maketitle

\section{Introduction}

An action of the torus $\T^k = (S^1)^k$ on a symplectic manifold
$(M,\omega)$ is a homomorphism
$$ \rho \colon \T^k \to \Sympl(M,\omega) $$
from the torus to the group of symplectomorphisms
such that the map $\rho^\dagger \colon \T^k \times M \to M$
given by $(a,m) \mapsto \rho(a)(m)$ is smooth.
A \emph{sub-torus} of $\Sympl(M,\omega)$ is the image of such a homomorphism.  
We restrict our attention to actions that are effective
(faithful: $\rho$ is one-to-one).

We say that two actions, $\rho_1 \colon \T^k \to \Sympl(M,\omega)$ 
and $\rho_2 \colon \T^k \to \Sympl(M,\omega)$, are \emph{equivalent} 
if there exists a symplectomorphism $\psi \colon M \to M$ 
and an automorphism $h \colon \T^k \to \T^k$ such that the diagram
$$\begin{CD}
 \T^k \times M @> {\rho_1}^\dagger >> M \\
@V (h,\psi) VV @V \psi VV \\
 \T^k \times M @> {\rho_2}^\dagger >> M 
\end{CD}$$
commutes.
Two (effective) torus actions on $(M,\omega)$ are equivalent 
if and only if their images are conjugate sub-tori of $\Sympl(M,\omega)$. 

A $\T^k$-action on $(M,\omega)$ with generating vector fields
$\xi_1, \ldots, \xi_k$ is \emph{Hamiltonian} if it has
a \emph{moment map}, that is, a map $ \Phi \colon M \to \R^k $
whose components satisfy Hamilton's equation
$d\Phi_j = -\iota(\xi_j) \omega$.
If $H^1(M;\R) = \{ 0 \}$ then every torus action on $(M,\omega)$ 
is Hamiltonian.  
An (effective) Hamiltonian torus action is \emph{toric} if the dimension
of the torus is half the dimension of the manifold;
if $M$ is compact and connected, the triple $(M,\omega,\Phi)$ is called 
a \emph{symplectic toric manifold}. 
In this paper we prove the following ``finiteness theorem":

\begin{Theorem} \labell{thm:finite}
Let $(M,\omega)$ be a four dimensional compact connected symplectic 
manifold.  Then the set of equivalence classes of toric actions 
on $(M,\omega)$ is finite.
\end{Theorem}

An important ingredient in the proof of Theorem \ref{thm:finite} is the
following properness property.
\begin{Lemma} \labell{sproper}
Let $(M,\omega)$ be a symplectic four-manifold that is obtained from
$\tCP$ by $k$ symplectic blow-ups.
Then the restriction of the period map $E \mapsto \omega(E)$ to the set
of classes in $H_2(M;\R)$ with self intersection $-1$ is proper.
\end{Lemma}
We prove Lemma \ref{sproper} and Theorem \ref{thm:finite} by ``soft''
equivariant, algebraic, and combinatorial techniques. There also exists
a properness result for a general compact symplectic four-manifold.
This is due to results of McDuff \cite[\S 3]{McD:topology},
\cite[Lemma~3.1]{McDuff-Structure}, using ``hard" holomorphic
techniques. In an appendix, we present these results and their proofs,
and apply them to (re)prove Theorem \ref{thm:finite}. This approach
exhibits the interplay between ``soft'' and ``hard''  techniques;
see \cite{kk} for another usage of holomorphic techniques in proving
results on toric actions.

Theorem \ref{thm:finite} has the following corollary.
\begin{Corollary} \labell{finite conj}
Let $(M,\omega)$ be a four dimensional compact connected symplectic manifold
with $H^1(M;\R) = \{ 0 \}$.   Then the
number of conjugacy classes of 2-tori in $\Sympl(M,\omega)$ is finite.
\end{Corollary}

The number of conjugacy classes can be greater than one; 
see Example \ref{hirzebruch}.

If $(M,\omega)$ admits a toric action then $M$ is simply connected
\cite[Corollary 32.2]{GSbook}, so $H^1(M;\R) = \{ 0 \}$.
We restrict our attention to symplectic manifolds $(M,\omega)$ 
such that $M$ is compact and connected and $H^1(M;\R) = \{ 0 \}$
and, as before, to torus actions that are effective.

Since the orbits of a Hamiltonian torus action are isotropic
(see \cite[Proposition III.2.12]{naudin}) and, on an open dense set, 
the (effective) action is free (see \cite[Corollary B.48]{GGK}), 
a toric action does not extend to a Hamiltonian action of a larger torus.  
Thus, the image of a toric action is a maximal torus in 
the symplectomorphism group $\Sympl(M,\omega)$.
If $\dim M = 4$, $H^1(M;\R) = \{ 0 \}$, and $\dim H^2(M;\R) \leq 3$, 
every Hamiltonian $S^1$-action on $(M,\omega)$ extends to a toric action 
\cite[Theorem~1]{karshon:max}, so in this case every maximal torus
in the symplectomorphism group comes from a toric action.
However, in general, there exist Hamiltonian $S^{1}$-actions that do not 
extend to toric actions, so there exist maximal tori that do not come
from toric actions.
In fact, some symplectic four-manifolds admit $S^1$-actions and do not
admit any toric actions.
For others, the symplectomorphism group contains some maximal tori 
that are two dimensional and some that are one dimensional.
See Example \ref{s1t}.
However, most symplectic four-manifolds do not admit any circle actions 
(in particular they do not admit any toric actions); 
see Example \ref{s1t}, and for many more examples see \cite{gompf}.

Maximal tori play a central role in the theory of finite dimensional 
Lie groups.  In a compact Lie group, every two maximal tori are conjugate. 
Our work can be viewed as an attempt to find aspects 
of finite dimensional Lie theory that carry over to groups of 
symplectomorphisms. See~\cite{Re}, \cite{McD-Kedra}, and \cite{bloch} 
for related ideas.
Homotopy theoretic properties of symplectomorphism groups have been studied
in \cite{abreu-mcduff,anjos-granja,buse,gromov,lalonde-pinsonnault,
seidel}; see the survey \cite{mcduff-survey}.

At this point we do not know whether a finiteness result
similar to Theorem \ref{thm:finite} also holds in higher dimensions
or for non-Hamiltonian symplectic group actions.  
The analogous result does \emph{not}
hold for contactomorphisms: Eugene Lerman \cite{lerman} 
constructed a compact contact $3$-manifold that admits 
infinitely many non-conjugate toric actions.
``Finiteness" also does not hold for the groups of all diffeomorphisms:
for example, $\Diff(S^2 \times S^2)$ contains infinitely many 
non-conjugate two dimensional maximal tori; see Example \ref{hirzebruch}.


\subsection*{Acknowledgements} We are grateful to Paul Biran, 
Dietmar Salamon, and Margaret Symington for helpful conversations.  
The project is partially funded by an NSERC Discovery grant, 
an NSERC LSI grant, an NSERC postdoctoral grant number BP-301203-2004, 
and a VIGRE grant DMS-9983190.

\section{Four dimensional symplectic toric manifolds and corner chopping} 
\labell{del}

The automorphism group of the torus $\T^n = \R^n / 2\pi \Z^n$
is the group $\GL(n,\Z)$ of $n \times n$ matrices with integer entries
and determinant equal to $1$ or $-1$.  
If a matrix $A$ belongs to $\GL(n,\Z)$ then so do 
its transpose $A^T$ and its inverse $A\Inv$.
We also consider the group $\AGL(n,\Z)$ of affine transformations
of $\R^n$ that have the form $x \mapsto Ax + \alpha$ 
with $A \in \GL(n,\Z)$ and $\alpha \in \R^n$.

A \emph{Delzant polytope} in $\R^n$ is a compact convex polytope $\Delta$ 
such that the edges emanating from each vertex 
can be generated by a $\Z$-basis of the lattice $\Z^n$.

\begin{Remark}
Other authors also use the following terms for polytopes with this property:  
non-singular, primitive, regular, smooth, unimodular, torsion free.
\end{Remark}

An important model for a toric action is $\C^n$ with the 
standard symplectic form, the standard $\T^n$-action given by rotations of the 
coordinates, and the moment map
\begin{equation} \labell{Cn mm}
 (z_1,\ldots,z_n) \mapsto \half (|z_1|^2,\ldots,|z_n|^2).
\end{equation}
The image of this moment map is the positive orthant, 
$$ \R_+^n = \{ (s_1,\ldots,s_n) \ | \ s_j \geq 0 \text{ for all $j$ } \}.$$

The moment map image of a symplectic toric manifold is a Delzant polytope.
A Delzant polytope can be obtained by gluing open subsets of $\R_+^n$ 
by means of elements of $\AGL(n,\Z)$.
Similarly, a symplectic toric manifold can be obtained by gluing 
open $\T^n$-invariant subsets of $\C^n$ by means of 
equivariant symplectomorphisms and reparametrizations of $\T^n$. 
In fact,

\begin{Lemma} \labell{local}
Let $(M,\omega)$ be a symplectic toric manifold with moment map 
$\Phi \colon M \to \Lie(\T^n)^* = \R^n$.
Let $x$ be a point in $\R^n$. Then there exists an open neighborhood
of $x$ whose moment map preimage is equivariantly symplectomorphic 
to an open subset of $\C^n$ with the standard $\T^n$-action
conjugated by some automorphism $\T^n \to \T^n$.
\end{Lemma} 

\begin{proof}
Because the moment map image is closed, 
the lemma is true if $\Phi\Inv(x)$ is empty.

Let $O$ be an orbit in $\Phi\Inv(x)$.
By the local normal form for Hamiltonian torus actions, 
some neighborhood $U$ of $O$ is equivariantly symplectomorphic
to an open subset of $\C^n$.  In particular, 
the level set $\Phi\Inv(x)$ is a discrete union of orbits.
Because the level sets are connected \cite{atiyah,GS:convexity}, 
$\Phi\Inv(x) = O$.
Because $\Phi$ is proper, the inverse image of a sufficiently small 
neighborhood of $x$ is contained in $U$.
\end{proof}

Delzant \cite{delzant} classified symplectic toric manifolds 
up to equivariant symplectomorphism.  Applying his theorem 
to different toric actions on a fixed manifold, we get the following result.

\begin{Proposition} \labell{prop:delzant}
Let $(M,\omega)$ be a $2n$ dimensional compact connected symplectic manifold.
\begin{enumerate}
\item \labell{d1}
For any Hamiltonian $\T^n$-action on $(M,\omega)$, its moment map image
is a Delzant polytope in~$\R^n$.
\item \labell{d2}
Two Hamiltonian $\T^n$-actions on $(M,\omega)$ are equivalent
if and only if their moment map images are $\AGL(n,\Z)$-congruent.
\end{enumerate}
\end{Proposition}

\begin{proof}
By \cite{atiyah,GS:convexity}, the moment map image is a convex polytope.
Lemma \ref{local} implies that this polytope is locally 
$\AGL(n,\Z)$-congruent to $\R_+^n$; thus, it is a Delzant polytope. 

To prove part \eqref{d2},
let $\rho_1 , \rho_2 \colon \R^n \to \Sympl(M,\omega)$
be toric actions with moment maps $\Phi_1 , \Phi_2 \colon M \to \R^n$
and moment map images $\Delta_1$ and $\Delta_2$.

Suppose that the actions $\rho_1$ and $\rho_2$ are equivalent,
i.e., that there exist $A \in \GL(n,\Z)$ and $\psi \in \Sympl(M,\omega)$
such that
\begin{equation} \labell{conj}
 \rho_2 (A \cdot \lambda) = \psi \circ \rho_1(\lambda) \circ \psi\Inv
\end{equation}
for all $\lambda \in \T^n$.
A moment map for the $\T^n$-action $\lambda \mapsto \rho_2(A \cdot \lambda)$
is $A^T \circ \Phi_2$.  A moment map for the $\T^n$-action 
$\lambda \mapsto \psi \circ \rho_1(\lambda) \circ \psi\Inv$ 
is $\Phi_1 \circ \psi\Inv$.
By \eqref{conj}, these are moment maps for the same action, 
so their difference is a constant: for some $\alpha \in \R^n$,
$$ \Phi_1 \circ \psi\Inv = A^T \circ \Phi_2 + \alpha . $$
Comparing the images of these maps, we get $\Delta_1 = A^T \Delta_2 + \alpha$.
So $\Delta_1$ and $\Delta_2$ are $\AGL(n,\Z)$-congruent.

Suppose that $\Delta_1$ and $\Delta_2$ are 
$\AGL(n,\Z)$-congruent.  Let $A \in \GL(n,\Z)$ and $\alpha \in \R^n$
be such that 
\begin{equation} \labell{transl}
\Delta_2 = A \Delta_1 + \alpha.
\end{equation}
The map $A \circ \Phi_1 + \alpha$ is a moment map for the action 
$\rho_1 \circ A^T$.  Its image, $A \Delta_1 + \alpha$, coincides 
with the image $\Delta_2$ of the moment map for the action $\rho_2$.  
By Delzant's theorem \cite{delzant}, toric actions with the same
moment map image are equivariantly symplectomorphic.  Applying this
to the actions $\rho_1 \circ A^T$ and $\rho_2$, we get that
there exists a symplectomorphism $\psi \colon M \to M$ such that 
$\psi \circ \rho_1 ( A^T \cdot \lambda) = \rho_2(\lambda) \circ \psi$
for all $\lambda \in \T^n$.
That is, the actions $\rho_1$ and $\rho_2$ are equivalent.
\end{proof}

\begin{noTitle} \labell{rational length}
The \emph{rational length} of an interval $d$ of rational slope in $\R^n$
is the unique number $\ell = \length(d)$ such that the interval 
is $\AGL(n,\Z)$-congruent to an interval of length $\ell$ 
on a coordinate axis.
In what follows, intervals are always measured by rational length.
\end{noTitle}

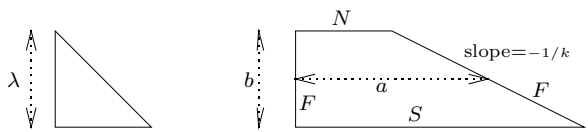
\begin{figure}[ht]
\setlength{\unitlength}{0.00083333in}
\begingroup\makeatletter\ifx\SetFigFont\undefined%
\gdef\SetFigFont#1#2#3#4#5{%
  \reset@font\fontsize{#1}{#2pt}%
  \fontfamily{#3}\fontseries{#4}\fontshape{#5}%
  \selectfont}%
\fi\endgroup%
{\renewcommand{\dashlinestretch}{30}
\begin{picture}(3612,639)(0,-10)

%
%
\path(300,612)(900,12)(300,12)(300,612)
\dottedline{45}(150,612)(150,12)
\path(120.000,132.000)(150.000,12.000)(180.000,132.000)
\path(180.000,492.000)(150.000,612.000)(120.000,492.000)
\put(0,260){\makebox(0,0)[lb]{\smash{{{\SetFigFont{12}{14.4}{\rmdefault}{\mddefault}{\updefault}$\scriptstyle{\lambda}$}}}}}
%
%
\path(1800,612)(2400,612)(3600,12)(1800,12)(1800,612)
%
%
\path(1920.000,342.000)(1800.000,312.000)(1920.000,282.000)
\dottedline{45}(1800,312)(3000,312)
\path(2880.000,282.000)(3000.000,312.000)(2880.000,342.000)
\put(2300,230){\makebox(0,0)[lb]{\smash{{{\SetFigFont{12}{14.4}{\rmdefault}{\mddefault}{\updefault}$\scriptstyle{a}$}}}}}
%
%
\dottedline{45}(1570,612)(1570,12)
\path(1600.000,492.000)(1570.000,612.000)(1560.000,492.000)
\path(1560.000,132.000)(1570.000,12.000)(1600.000,132.000)
\put(1480,260){\makebox(0,0)[lb]{\smash{{{\SetFigFont{12}{14.4}{\rmdefault}{\mddefault}{\updefault}$\scriptstyle{b}$}}}}}
\put(1820,120){\makebox(0,0)[lb]{\smash{{{\SetFigFont{12}{14.4}{\rmdefault}{\mddefault}{\updefault}$\scriptstyle{F}$}}}}}
%
%
\put(2850,450){\makebox(0,0)[lb]{\smash{{{\SetFigFont{12}{14.4}{\rmdefault}{\mddefault}{\updefault}\scriptsize{slope=}$\scriptstyle{-1/k}$}}}}}
%
%
\put(2500,50){\makebox(0,0)[lb]{\smash{{{\SetFigFont{12}{14.4}{\rmdefault}{\mddefault}{\updefault}$\scriptstyle{S}$}}}}}
\put(2025,650){\makebox(0,0)[lb]{\smash{{{\SetFigFont{12}{14.4}{\rmdefault}{\mddefault}{\updefault}$\scriptstyle{N}$}}}}}
\put(3275,200){\makebox(0,0)[lb]{\smash{{{\SetFigFont{12}{14.4}{\rmdefault}{\mddefault}{\updefault}$\scriptstyle{F}$}}}}}
\end{picture}
}
\caption{A Delzant triangle, $\Delta_\lambda$, 
         and a Hirzebruch trapezoid, $H_{a,b,k}$}
\labell{fig:triangle}
\end{figure}

\begin{Example} \labell{ex:Delzant}
Figure \ref{fig:triangle} shows examples of Delzant polygons 
with three and four edges.  
On the left there is a \emph{Delzant triangle},
$$ \Delta_\lambda = \{ (x_1,x_2) \ | \ 
   x_1 \geq 0 , x_2 \geq 0 , x_1 + x_2 \leq \lambda \}. $$
This is the moment map image of the standard toric action on $\CP^2$, 
with the Fubini-Study symplectic form normalized so that 
the symplectic area of $\CP^1 \subset \CP^2$ is $2 \pi \lambda$.
The rational lengths of all its edges is $\lambda$.

On the right there is a \emph{Hirzebruch trapezoid},
$$ H_{a,b,k} = \{ (x_1,x_2) \ | \ 
   -\frac{b}{2} \leq x_2 \leq \frac{b}{2} , 0 \leq x_1 \leq a - kx_2 \} .$$
$b$ is the height of the trapezoid, $a$ is its average width, 
and $k$ is a non-negative integer such that the right edge has slope $-1/k$ 
or is vertical if $k=0$.  
We assume that $a \geq b$ and that $a - k \frac{b}{2} >  0$.
This trapezoid is a moment map image of a standard toric action 
on a Hirzebruch surface.
The rational lengths of its left and right edges are $b$;
the rational lengths of its top and bottom edges are $a \pm kb/2$.
\eoe
\end{Example}

\begin{Example} \labell{hirzebruch}
Consider the manifold $S^2 \times S^2$ with the product symplectic form
$\omega_{a,b} = a \tau \oplus b \tau$, where $\tau$ is the rotation 
invariant area form on $S^2$ with total area $2\pi$.
For any toric action on $(S^2 \times S^2 , \omega_{a,b})$, 
its moment map image is $\AGL(2,\Z)$-congruent to the 
Hirzebruch trapezoid $H_{a,b,k}$ for some non-negative even integer $k$.
Given $a$ and $b$, there exist only finitely many such $k$'s for which
$a - k \frac{b}{2} > 0$.  Hence, up to conjugation,
there are only finitely many 2-tori in the symplectomorphism group
of $(S^2 \times S^2, \omega_{a,b})$.
See \cite{karshon:max}.

Different non-negative even integers $k$ 
correspond to $2$-tori in the symplectomorphism group
$\Sympl(S^2 \times S^2, \omega_{a,b})$
that are not conjugate.  
$2$-tori corresponding to different integers $k$ are even not conjugate
in the group
$\Diff(S^2 \times S^2)$ of \emph{diffeomorphisms} of $S^2 \times S^2$:
the value of $k$ can be recovered from the fact that the two dimensional 
orbit type strata have self intersection $0$, $k$, and $-k$.
These tori are maximal in $\Diff(S^2 \times S^2)$: because an
effective compact group action on a connected manifold is also
effective on invariant open subsets, and by local linearization, a
torus that acts effectively on a manifold with a fixed point is at
most of half the dimension of the manifold.  
Hence, $\Diff(S^2 \times S^2)$ contains \emph{infinitely many} 
conjugacy classes of maximal tori.
\eoe
\end{Example}

\begin{noTitle} \labell{combinatorial self intersection}
Let  $d$, $d'$, and $d''$ be three consecutive edges in a Delzant polygon, 
ordered counterclockwise.  Let  $u$, $u'$, and $u''$ be their outward 
normal vectors, normalized so that they are primitive lattice elements. 
Then each of $(u,u')$ and $(u',u'')$ is an oriented $\Z$-basis 
for $\Z^2$.  It follows that there exists an integer  $k$ such that  
$u + u'' = k u'$. 
We define the \emph{combinatorial self intersection number} of $d'$ 
to be $-k$. We also say that the 
\emph{combinatorial intersection number} of two different edges  
is $1$ if they are adjacent and $0$ if not.
\end{noTitle}

\begin{figure}[ht]
\setlength{\unitlength}{0.0006in}
\begingroup\makeatletter\ifx\SetFigFont\undefined%
\gdef\SetFigFont#1#2#3#4#5{%
  \reset@font\fontsize{#1}{#2pt}%
  \fontfamily{#3}\fontseries{#4}\fontshape{#5}%
  \selectfont}%
\fi\endgroup%
{\renewcommand{\dashlinestretch}{30}
\begin{picture}(1824,2810)(0,-10)
\put(612,312){\blacken\ellipse{70}{70}}
\put(912,312){\blacken\ellipse{70}{70}}
\put(1212,312){\blacken\ellipse{70}{70}}
\put(1512,312){\blacken\ellipse{70}{70}}
\put(312,612){\blacken\ellipse{70}{70}}
\put(612,612){\blacken\ellipse{70}{70}}
\put(912,612){\blacken\ellipse{70}{70}}
\put(1212,612){\blacken\ellipse{70}{70}}
\put(1512,612){\blacken\ellipse{70}{70}}
\put(312,912){\blacken\ellipse{70}{70}}
\put(612,912){\blacken\ellipse{70}{70}}
\put(912,912){\blacken\ellipse{70}{70}}
\put(1212,912){\blacken\ellipse{70}{70}}
\put(1512,912){\blacken\ellipse{70}{70}}
\put(312,1212){\blacken\ellipse{70}{70}}
\put(612,1212){\blacken\ellipse{70}{70}}
\put(912,1212){\blacken\ellipse{70}{70}}
\put(1212,1212){\blacken\ellipse{70}{70}}
\put(1512,1212){\blacken\ellipse{70}{70}}
\put(312,1512){\blacken\ellipse{70}{70}}
\put(612,1512){\blacken\ellipse{70}{70}}
\put(912,1512){\blacken\ellipse{70}{70}}
\put(1212,1512){\blacken\ellipse{70}{70}}
\put(1512,1512){\blacken\ellipse{70}{70}}
\put(312,1812){\blacken\ellipse{70}{70}}
\put(612,1812){\blacken\ellipse{70}{70}}
\put(912,1812){\blacken\ellipse{70}{70}}
\put(312,312){\blacken\ellipse{70}{70}}
\put(1212,1812){\blacken\ellipse{70}{70}}
\put(1512,1812){\blacken\ellipse{70}{70}}
\put(312,2112){\blacken\ellipse{70}{70}}
\put(612,2112){\blacken\ellipse{70}{70}}
\put(912,2112){\blacken\ellipse{70}{70}}
\put(1212,2112){\blacken\ellipse{70}{70}}
\put(1512,2112){\blacken\ellipse{70}{70}}
\put(312,2412){\blacken\ellipse{70}{70}}
\put(612,2412){\blacken\ellipse{70}{70}}
\put(912,2412){\blacken\ellipse{70}{70}}
\put(1212,2412){\blacken\ellipse{70}{70}}
\put(1512,2412){\blacken\ellipse{70}{70}}
\put(312,2712){\blacken\ellipse{70}{70}}
\put(612,2712){\blacken\ellipse{70}{70}}
\put(912,2712){\blacken\ellipse{70}{70}}
\put(1212,2712){\blacken\ellipse{70}{70}}
\put(1512,2712){\blacken\ellipse{70}{70}}
\path(312,2712)(912,2112)(1512,912)
	(1512,312)(312,312)(312,2712)
\path(912,312)(912,12)
\path(882.000,132.000)(912.000,12.000)(942.000,132.000)
\path(1512,612)(1812,612)
\path(1692.000,582.000)(1812.000,612.000)(1692.000,642.000)
\path(1212,1512)(1812,1812)
\path(1718.085,1731.502)(1812.000,1812.000)(1691.252,1785.167)
\path(612,2412)(912,2712)
\path(848.360,2605.934)(912.000,2712.000)(805.934,2648.360)
\path(312,1512)(12,1512)
\path(132.000,1542.000)(12.000,1512.000)(132.000,1482.000)

\put(662,362){\makebox(0,0)[lb]{\smash{{{\SetFigFont{12}{14.4}{\rmdefault}{\mddefault}{\updefault}$\sss d_1$}}}}}
\put(1307,467){\makebox(0,0)[lb]{\smash{{{\SetFigFont{12}{14.4}{\rmdefault}{\mddefault}{\updefault}$\sss d_2$}}}}}
\put(937,1627){\makebox(0,0)[lb]{\smash{{{\SetFigFont{12}{14.4}{\rmdefault}{\mddefault}{\updefault}$\sss d_3$}}}}}
\put(582,2187){\makebox(0,0)[lb]{\smash{{{\SetFigFont{12}{14.4}{\rmdefault}{\mddefault}{\updefault}$\sss d_4$}}}}}
\put(362,1087){\makebox(0,0)[lb]{\smash{{{\SetFigFont{12}{14.4}{\rmdefault}{\mddefault}{\updefault}$\sss d_5$}}}}}
\end{picture}

}
\caption{}
\labell{fig:polygon}
\end{figure}

\begin{Example} \labell{polygon}
The edges of the polygon in figure \ref{fig:polygon}
have the following rational lengths and combinatorial self intersections:
\begin{center}
\begin{tabular}{c|c|c|}
& rational length & combinatorial self intersection \\
\hline
$d_1$ & 4 & \phantom{-}0 \\
\hline
$d_2$ & 2 & -2 \\
\hline
$d_3$ & 2 & -1 \\
\hline
$d_4$ & 2 & -1 \\
\hline
$d_5$ & 8 & \phantom{-}1  \\
\hline
\end{tabular}
\end{center}
\eoe
\end{Example}

\begin{noTitle} \labell{calJ}
An almost complex structure on a $2n$ dimensional manifold $M$ is an
automorphism of the tangent bundle, $J \colon TM \to TM$,
such that $J^2 = -\id$.  It is \emph{tamed} by a symplectic form $\omega$ 
if $\omega(v,Jv) > 0$ for all $v \neq 0$. The \emph{first Chern class} 
of the symplectic manifold $(M,\omega)$ is defined to be the first Chern 
class of the complex vector bundle $(TM,J)$ and is denoted $c_1(TM)$.
This class is independent of the choice of tamed almost complex structure $J$
\cite[\S 2.6]{MS:intro}.
\end{noTitle}

\begin{Lemma} \labell{perimeter and area}
Let $(M,\omega)$ be a compact connected symplectic four-manifold.  
Let $\Phi \colon M \to \R^n$  be the moment map for a toric action,
and let $\Delta$ be its image.
\begin{enumerate}
\item \labell{i2}
Let $d$ be an edge of $\Delta$ of rational length $\ell$.  
Then its preimage, $\Phi\Inv(d)$, is a symplectically embedded 2-sphere 
in $M$ of symplectic area 
$$ \int_{\Phi\Inv(d)} \omega = 2 \pi \ell.$$


\item \labell{i2'} 
Let $d$ and $d'$ be edges of $\Delta$. Then their combinatorial 
intersection number (see \S\ref{combinatorial self intersection})
is equal to the intersection number of their 
preimages in $M$. 
\item \labell{i1}
The preimages of the edges of $\Delta$ generate the second homology 
group of $M$. 
The number of vertices of $\Delta$ is equal to $\dim H_2(M) + 2$.
\item \labell{i3}
The perimeter of $\Delta$ (see \S\ref{rational length}) is
$$ \perimeter(\Delta) = \frac{1}{2\pi} \int_M \omega \wedge c_1(TM) .$$
\item \labell{i4}
The area of $\Delta$ is
$$ \frac{1}{(2\pi)^2} \int_M \frac{1}{2!} \omega \wedge \omega .$$
\end{enumerate}
\end{Lemma}

\begin{proof}
By Lemma \ref{local}, the preimage of a sufficiently small subinterval
$c$ of the edge $d$ is a disk or an annulus with area $2\pi \length(c)$.
Part (\ref{i2}) follows.

When $d \neq d'$, part (\ref{i2'}) follows from Lemma \ref{local}. 
When $d=d'$, let $d_-$ and $d_+$ be the edges that are adjacent to $d$.
Let $u_-$, $u$, and $u_+$ be the primitive outward normal vectors
to the edges $d_-$, $d$, and $d_+$. Consider the sub-circle 
$i \colon S^1 \to \T^2$ determined by $u_- \in \Lie(\T^2)$.
Then $i(S^1)$ acts on the sphere $\Phi\Inv(d)$ by rotations with speed 1,
it rotates the normal to the sphere at the north pole with speed $k$, 
where $-k$ is the combinatorial self intersection of $d$,
and it fixes the normal to the sphere at the south pole.
By examining the transition function for a trivializations
of the normal bundle over complements of the north and south poles,
we see that the Euler number of this normal bundle is $-k$.

By, e.g., Lemma \ref{local},
a generic component of the moment map is a Morse function
with even indices whose critical points are the preimages of the 
vertices of the polytope.  For a critical point of index two,
the descending gradient manifold with respect to an invariant metric
is the preimage of an edge of $\Delta$ minus its bottom vertex.
Part (\ref{i1}) follows from this using Morse theory.

Let $C_1,\ldots,C_N$ be the preimages in $M$
of the edges of $\Delta$.  Then the homology class
$\sum_{i=1}^N [C_i] \in H_2(M;\R)$ is the Poincar\'e dual
to the first Chern class $c_1(TM)$.
This follows from the fact that the $C_i$'s generate the second homology
group and since
\begin{align*}
[C_i]  \cdot \sum_{j=1}^N [C_j] &= [C_i] \cdot [C_i] + 2  
      \qquad \text{ by part \eqref{i2'} }\\
 &= {c_1(TM)([C_i])}  
      \qquad \text{ since $C_i$ is an embedded 2-sphere in a four-manifold.}
 \end{align*}
This and part \eqref{i2} imply part \eqref{i3}.

Part \eqref{i4} follows from Lemma~\ref{local}.
It is a special case of the Duistermaat-Heckman theorem.
\end{proof}

\begin{noTitle} \labell{corner chopping}
Let $\Delta$ be a Delzant polytope in $\R^n$,
let $v$ be a vertex of $\Delta$,
and let $\delta > 0$ be smaller than the rational lengths
of the edges emanating from $v$. The edges of $\Delta$
emanating from $v$ have the form 
$\{ v + s \alpha_j \ | \ 0 \leq s \leq \ell_j \}$
where the vectors $\alpha_1, \ldots \alpha_n$ generate the lattice $\Z^n$
and $\delta < \ell_j$ for all $j$.
The \emph{corner chopping of size $\delta$} of $\Delta$ at $v$
is the polytope $\tDelta$ obtained from $\Delta$ by intersecting
with the half-space
$$ \{ \ v + s_1 \alpha_1 + \ldots + s_n \alpha_n \ \ | \  \ \ 
      s_1 + \ldots + s_n \geq \delta \} . $$
See, e.g., the chopping of the top right corner in Figure \ref{fig:blowup}.
The resulting polytope $\tDelta$ is again a Delzant polytope.
The ``corner chopping" operation commutes with $\AGL(n,\Z)$-congruence:
if $\tDelta$ is obtained from $\Delta$ by a corner chopping
of size $\delta > 0$ at a vertex $v \in \Delta$ then, for any
$g \in \AGL(n,\Z)$, the polytope $g(\tDelta)$ is obtained
from the polytope $g(\Delta)$ by a corner chopping of size $\delta$
at the vertex~$g(v)$.  
\end{noTitle}

\begin{Remark} \labell{blowup}
If $M$ is a toric manifold with moment map image $\Delta$,
and $\tDelta$ can be obtained from $\Delta$ by a corner chopping,
then $\tDelta$ is the moment map image of a toric manifold
obtained from $M$ by an equivariant symplectic blow-up.
See~\cite[Section 3]{kk}.
\end{Remark}

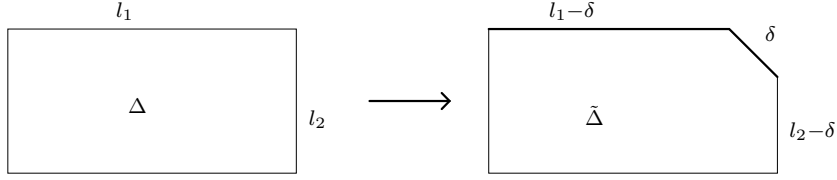
\begin{figure}[ht]
\setlength{\unitlength}{0.00083333in}
\begingroup\makeatletter\ifx\SetFigFont\undefined%
\gdef\SetFigFont#1#2#3#4#5{%
  \reset@font\fontsize{#1}{#2pt}%
  \fontfamily{#3}\fontseries{#4}\fontshape{#5}%
  \selectfont}%
\fi\endgroup%
{\renewcommand{\dashlinestretch}{30}
\begin{picture}(5211,1137)(0,-10)
\thinlines
\path(3044,912)(4544,912)(4844,612)
	(4844,12)(3044,12)(3044,912)
\thicklines
\path(3044,912)(4544,912)(4844,612)
\put(719,987){\makebox(0,0)[lb]{\smash{{{\SetFigFont{12}{14.4}{\rmdefault}{\mddefault}{\updefault}$\sss l_1$}}}}}
\path(2300,462)(2800,462)
\path(2750,512)(2800,462)(2750,412)
\thinlines
\path(44,912)(1844,912)(1844,12)
	(44,12)(44,912)
\put(1919,312){\makebox(0,0)[lb]{\smash{{{\SetFigFont{12}{14.4}{\rmdefault}{\mddefault}{\updefault}$\sss l_2$}}}}}
\put(4919,237){\makebox(0,0)[lb]{\smash{{{\SetFigFont{12}{14.4}{\rmdefault}{\mddefault}{\updefault}$\sss l_2-\delta$}}}}}
\put(794,387){\makebox(0,0)[lb]{\smash{{{\SetFigFont{12}{14.4}{\rmdefault}{\mddefault}{\updefault}$\sss \Delta$}}}}}
\put(3644,312){\makebox(0,0)[lb]{\smash{{{\SetFigFont{12}{14.4}{\rmdefault}{\mddefault}{\updefault}$\sss \Tilde{\Delta}$}}}}}
\put(3419,987){\makebox(0,0)[lb]{\smash{{{\SetFigFont{12}{14.4}{\rmdefault}{\mddefault}{\updefault}$\sss l_1-\delta$}}}}}
\put(4769,837){\makebox(0,0)[lb]{\smash{{{\SetFigFont{12}{14.4}{\rmdefault}{\mddefault}{\updefault}$\sss \delta$}}}}}
\end{picture}
}
\caption{Corner chopping} 
\labell{fig:blowup}
\end{figure}

The following lemma is an easy consequence of the definitions
of corner chopping and of combinatorial self intersection.

\begin{Lemma} \labell{intersection}
Let  $d$  and  $d'$  be consecutive edges in a Delzant polygon $\Delta$.
Let  $\tilde{\Delta}$  be obtained from  $\Delta$  by a corner chopping
at the vertex between  $d$  and  $d'$.  Let  $e$  be the resulting new edge.
Let  $\tilde{d}$  and  $\tilde{d'}$ be the edges of  $\tilde{\Delta}$  with
the same outward normal vectors as those of  $d$  and  $d'$.  Then 
\begin{enumerate}
\item the combinatorial self intersection of $e$ is $-1$,
\item the combinatorial self intersection of  $\tilde{d}$  is one less
than that of  $d$,
\item the combinatorial self intersection of  $\tilde{d'}$  is one less
than that of  $d'$,  
\item for any of the other edges of $\tilde{\Delta}$, 
its combinatorial self intersection is the same as that 
of the edge of $\Delta$ with
the same outward normal vector.
\end{enumerate}
\end{Lemma}

\begin{Remark}
Let $(M,\omega,\Phi)$ be a toric manifold with moment map image $\Delta$.
Lemma \ref{intersection} is the combinatorial counterpart of the following
geometric facts.   Let $\tM$ be the toric manifold obtained from $M$
by equivariant symplectic blow-up centered at a point $p$. 
The exceptional divisor in $\tM$ has self intersection $-1$; a 2-sphere 
through $p$ in $M$ with self intersection $\ell$
gives rise to a 2-sphere in $\tM$ with self intersection $\ell - 1$.
\end{Remark}

\begin{Remark} \labell{s2triv}
Recall that for any toric action on $(S^2 \times S^2 , \omega_{a,b})$, 
its moment map image is $\AGL(2,\Z)$-congruent to the 
Hirzebruch trapezoid $H_{a,b,k}$ for some non-negative even integer $k$ 
(see Example \ref{hirzebruch});
the combinatorial self intersections of the edges of $H_{a,b,k}$ are, 
from the right edge in counterclockwise direction, $(0,-k,0,k)$.
This fact and part (1) of Lemma \ref{intersection} show that 
for any toric action on $(S^2 \times S^2, \omega_{a,b})$, 
the moment map image \emph{cannot} be obtained by a corner chopping 
from another Delzant polygon. 
\end{Remark}

In $\R^2$, all Delzant polygons can be obtained by a simple recursive recipe:

\begin{Lemma} \labell{fulton}
\begin{enumerate}
\item \labell{f1}
Let $\Delta$ be a Delzant polygon with three edges. Then there exists 
a unique $\lambda > 0$ such that $\Delta$ is $\AGL(2,\Z)$-congruent
to the Delzant triangle $\Delta_\lambda$.
(See Example \ref{ex:Delzant}.)
\item \labell{f2}
Let $\Delta$ be a Delzant polygon with four or more edges.
Let $s$ be the non-negative integer such that the number of edges
is $4+s$.  Then there exist positive numbers $a \geq b > 0$, 
an integer $0 \leq k \leq 2a/b$, and positive numbers
$\delta_1, \ldots,\delta_s$, such that $\Delta$ is $\AGL(2,\Z)$-congruent 
to a Delzant polygon that is obtained from the Hirzebruch trapezoid 
$H_{a,b,k}$ 
(see Example \ref{ex:Delzant})
by a sequence of corner choppings of sizes $\delta_1, \ldots, \delta_s$.
\item \labell{f3}
Moreover, let $\Delta$ be a Delzant polygon with five or more edges.
Then $\Delta$ is $\AGL(2,\Z)$-congruent to a Delzant polygon that is
obtained from a Hirzebruch trapezoid $H_{a,b,k}$ with $k$ \emph{odd}
by a sequence of corner choppings. 
\end{enumerate}
\end{Lemma}

\begin{figure}[ht]
\setlength{\unitlength}{0.0006in}
\begingroup\makeatletter\ifx\SetFigFont\undefined%
\gdef\SetFigFont#1#2#3#4#5{%
  \reset@font\fontsize{#1}{#2pt}%
  \fontfamily{#3}\fontseries{#4}\fontshape{#5}%
  \selectfont}%
\fi\endgroup%
{\renewcommand{\dashlinestretch}{30}
\begin{picture}(3949,989)(0,-10)
\path(612,12)(612,612)(962,962)
	(2012,962)(3912,12)(612,12)
\dottedline{60}(962,962)(612,962)(612,612)
	(12,12)(612,12)
\end{picture}
}
\caption{}
\labell{even-odd}
\end{figure}

\begin{proof}
For parts \eqref{f1} and \eqref{f2}, see 
\cite[Section~2.5 and Notes to Chapter 2]{fulton}.
For part \eqref{f3}, let $\Delta$ be a Delzant polygon with five edges.

By part \eqref{f2}, up to $\AGL(2,\Z)$-congruence, $\Delta$ is obtained
from a Hirzebruch trapezoid $H_{a,b,k}$ by one corner chopping.
Suppose, without loss of generality, that the corner chopping is at the
left hand side.
The edges on the left and right hand sides of $H_{a,b,k}$ have
combinatorial self-intersection $0$,
hence, by Lemma \ref{intersection}, 
in $\Delta$ both the new edge created in the corner chopping and the
remainder $\tilde{f_l}$ of the left edge have self intersection $-1$.
So $\Delta$ can be obtained by another corner chopping, (in which
$\tilde{f_l}$ is the new edge), from a Delzant polygon that is
$\AGL(2,\Z)$-congruent to a Hirzebruch trapezoid with integer $k'$,
such that $k'=k+1$ or $k'=k-1$, depending on whether the corner
chopping was at the top or at the bottom left corner.
See Figure \ref{even-odd}.
\end{proof}


\begin{Corollary} \labell{enough}
Let $(M,\omega)$ be a compact connected symplectic four-manifold that
admits a toric action.  Then $(M,\omega)$ is symplectomorphic either to
$(S^2 \times S^2, \omega_{a,b})$ or to a manifold that is obtained
from  $\CP^2$  with a multiple of the Fubini-Study form by a sequence
of  $H_2(M;\R)-1$  symplectic blow-ups.
\end{Corollary}

\begin{proof}
Let $\Delta$ be the moment map image for a toric action on
$(M,\omega)$.  By part \eqref{d1} of Proposition~\ref{prop:delzant},
$\Delta$ is a Delzant polygon.  By part \eqref{i1} of Lemma
\ref{perimeter and area}, $\Delta$ has $\dim H_2(M;\R)+2$ edges.

If $\dim H_2(M;\R)=1$, then $\Delta$ has $3$ edges.
By part \eqref{f1} of Lemma \ref{fulton}, $\Delta$ is
$\AGL(2,\Z)$-congruent to a standard Delzant triangle,
$\Delta_\lambda$. Hence, by the first part of Example \ref{ex:Delzant},
$\Delta$ is the moment map image of a toric action on $\CP^2$ with a
multiple of the Fubini-Study form.

If $\dim H_2(M;\R)=2$, then $\Delta$ has $4$ edges.
By part \eqref{f2} of Lemma \ref{fulton}, 
$\Delta$ is $\AGL(2,\Z)$-congruent to a Hirzebruch
trapezoid $H_{a,b,k}$. Hence, by the second part of Example
\ref{ex:Delzant}, $\Delta$ is the moment map image of a toric action on
a Hirzebruch surface.  If $k=0$, the Hirzebruch surface 
is $(S^2 \times S^2, \omega_{a,b})$; if $k=1$, then $H_{a,b,k}$ can be
obtained by one corner chopping from a Delzant triangle, i.e., the
corresponding Hirzebruch surface is obtained from $\tCP$ by an
equivariant symplectic blow-up (see Remark \ref{blowup}).  
If $ k \geq 2$, then, by \cite[Lemma 3]{karshon:max}, the corresponding
Hirzebruch surface is symplectomorphic to the Hirzebruch surface which
corresponds to the integer $k-2$ with the same parameters $a$ and $b$,
so, by recursion, it is symplectomorphic 
either to $(S^2 \times S^2, \omega_{a,b})$ or to a symplectic manifold
obtained from $\tCP$ by a symplectic blow-up.

If $\dim H_2(M;\R) \geq 3$, then $\Delta$ has five or more edges.  By
part \eqref{f3} of Lemma \ref{fulton}, up to $\AGL(2,\Z)$-congruence,
$\Delta$ can be obtained by corner chopping from a Hirzebruch trapezoid
$H_{a,b,k}$ with $k$ odd.  By the argument above, 
every such trapezoid is a moment map image for some
toric action on a symplectic blow-up of $\CP^2$.

By Delzant's uniqueness theorem \cite{delzant}, two symplectic toric
manifolds are equivariantly symplectomorphic if and only if their
moment map images differ by a translation. This completes our proof.
\end{proof}

\begin{Lemma} \labell{bound}
Let 
$$ \Delta_0,\Delta_1,\ldots,\Delta_s $$
be a sequence of Delzant polygons such that, for each $i$,
the polygon $\Delta_i$ is obtained from the polygon $\Delta_{i-1}$
by a corner chopping of size $\delta_i$. Then
\begin{enumerate}
\item \labell{j1}
$ \area(\Delta_s) = \area(\Delta_0) 
                  - \half \delta_1^2 - \ldots - \half \delta_s^2 . $
\item \labell{j2}
$ \perimeter(\Delta_s) = \perimeter(\Delta_0) 
                       - \delta_1 - \ldots - \delta_s. $
\item \labell{j3}
For each $i$, \ 
$ \delta_i < 2^s \perimeter(\Delta_s).$
\item \labell{j4}
The length of each edge of each $\Delta_i$
is a linear combination, with non-negative integer coefficients, 
of the lengths of the edges of $\Delta_s$.
\end{enumerate}
\end{Lemma}

\begin{proof}
Parts \eqref{j1} and \eqref{j2} follow immediately from the definition
of corner chopping in \S\ref{corner chopping}.

Because a corner chopping of size $\delta_i$ can only be performed 
if there are two consecutive edges in $\Delta_{i-1}$ 
of lengths greater than $\delta_i$,
$$ \perimeter(\Delta_{i-1}) > 2 \delta_i .$$
Because a corner chopping of size $\delta_i$ creates a new edge 
of size $\delta_i$, shortens by $\delta_i$ the lengths of two edges, 
and the other edges do not change,
$$ \perimeter(\Delta_i) = \perimeter(\Delta_{i-1}) - \delta_i .$$
From these two equalities it follows that
$$ \half \perimeter(\Delta_{i-1}) < \perimeter(\Delta_i) .$$
By induction, we get 
$$ \perimeter (\Delta_i) < 2^{s-i} \perimeter(\Delta_s).$$
Because $\Delta_i$ contains an edge of size $\delta_i$,
$$ \delta_i < \perimeter(\Delta_i) .$$
The last two equalities imply part \eqref{j3}.

Let $d$ be an edge of $\Delta_{i-1}$, let $\tilde{d}$ be the edge 
of $\Delta_i$ with the same outward normal vector, and let $e$ be 
the edge of $\Delta_i$ that was created at the $i$-th corner chopping.
If $d$ does not touch the vertex of $\Delta_{i-1}$ at which we chop,
then $\length({d}) = \length(\tilde{d})$.
If $d$ touches the vertex of $\Delta_{i-1}$ at which we chop, then
$\length(d) = \length(\tilde{d}) + \length(e)$.
Part \eqref{j4} of the lemma follows by induction.
\end{proof}

\begin{noTitle} \labell{Z of edges}
For any Delzant polygon $\Delta$, consider the free Abelian group
generated by its edges:
\begin{equation} \labell{edge group}
 \Z[\text{edges of }\Delta] .
\end{equation}

The ``combinatorial intersection pairing" on this group
is the $\Z$-valued bilinear pairing whose restriction to the generators
is given by their ``combinatorial intersection number".
The ``length functional"
$$ \Z[\text{edges of }\Delta] \to \R $$
is the homomorphism that associates to each basis element its rational 
length.
If $\Delta_{i+1}$ is obtained from $\Delta_i$ by a corner chopping,
we consider the injective homomorphism 
\begin{equation} \labell{injection}
 \Z[\text{edges of }\Delta_{i}] \hookrightarrow
 \Z[\text{edges of }\Delta_{i+1}] 
\end{equation}
whose restriction to the generators is defined in the following way.
If $d$ is an edge of $\Delta_i$ that does not touch the corner
that was chopped, then $d$ is mapped to the edge of $\Delta_{i+1}$
with the same outward normal vector.  If $d$ is an edge of $\Delta_i$
that touches the corner that was chopped, then $d$ is mapped
to $\tilde{d}+e$ where $e$ is the new edge of $\Delta_{i+1}$, 
created in the chopping, and $\tilde{d}$ is the edge of $\Delta_{i+1}$
with the same outward normal vector as $d$.
Lemma \ref{intersection} implies that the homomorphism \eqref{injection} 
respects the combinatorial intersection pairings.
The definition of corner chopping in \S\ref{corner chopping} 
implies that the homomorphism \eqref{injection} respects 
the length homomorphisms.
\end{noTitle}

\begin{Lemma} \labell{homomorphism}
If $M$ is a symplectic toric manifold with moment map image $\Delta$,
there exists a natural homomorphism from $\Z[\text{edges of }\Delta]$
onto $H_2(M;\Z)$: 
send a generator $d$ to the homology class of the 
two-sphere $\Phi\Inv(d) \subset M$.
This homomorphism carries the combinatorial intersection pairing
on $\Z[\text{edges of }\Delta]$ to the ordinary intersection pairing 
in homology, and it carries the length functional to the functional
on homology given by pairing with the symplectic form divided by $2\pi$.
\end{Lemma}

\begin{proof}
The lemma follows from parts (1)--(4) of Lemma \ref{perimeter and area}.
\end{proof}

\begin{Remark} \labell{natural}
Suppose that $\Delta_{i+1}$ is obtained from $\Delta_i$
by a corner chopping.  Let $M_i$ and $M_{i+1}$ 
be toric manifolds with moment map images $\Delta_{i}$.
and $\Delta_{i+1}$. Then $M_{i+1}$ can be obtained from $M_i$ 
by equivariant symplectic blow-up. 
The homomorphism \eqref{injection} induces the 
natural inclusion map of $H_2(M_i)$ into $H_2(M_{i+1})$.
\end{Remark} 

\section{Soft proofs of properness and finiteness} 
\labell{corner}

We now state the key lemma used in the proof of our main theorem:

\begin{Lemma} \labell{finite tuples}
Let $(M,\omega)$ be a closed connected symplectic four-manifold.
Then there exists a finite set $\calD(M,\omega)$ of positive numbers
such that for any polygon $\Delta$,
if $\Delta$ is the moment map image of a toric action on $(M,\omega)$, 
and if $\Delta$ can be obtained from a Delzant polygon $\Delta'$ 
by a sequence of corner choppings of sizes $\delta_1,\ldots,\delta_s$,
then $\delta_i \in \calD(M,\omega)$ for all $i$.
\end{Lemma}

\begin{noTitle} \labell{eqbound}
The proof of Lemma \ref{finite tuples} will use the following
observation. Let $\Delta$ be a polygon that satisfies the assumptions
of Lemma \ref{finite tuples}.
By part \ref{i1} of Lemma \ref{perimeter and area},
the number of vertices of $\Delta$ is $\dim H_2(M;\R) + 2$, so 
\begin{equation} \labell{bound on s}
s < \dim H_2(M;\R).
\end{equation}
By part \eqref{j3} of Lemma \ref{bound}, by part \eqref{i3} 
of Lemma~\ref{perimeter and area}, and by \eqref{bound on s},
$\delta_i$ is smaller than 
\begin{equation} \labell{bound on delta i}
2^{\dim H_2(M;\R)} \cdot \frac{1}{2\pi} \int_M \omega \wedge c_1(TM)  
\end{equation}
for all $i$.
\end{noTitle} 

When the cohomology class $[\omega]$ is rational,
that is, when there exists a positive number $h$ such that $\int_A \omega$ 
is an integer multiple of $h$ for all $A \in H_2(M;\Z)$,
one can prove Lemma \ref{finite tuples} as follows.
After replacing $\omega$ by $\frac{2\pi}{h} \omega$ 
(which does not affect the symplectomorphism group)
we may assume that $\frac{1}{2\pi} \int_A \omega \in \Z$
for all $A \in H_2(M;\Z)$. 

By part \eqref{i2} of Lemma \ref{perimeter and area}, 
the lengths of all the edges of $\Delta$ are integers.
By part \eqref{j4} of Lemma \ref{bound}, 
this implies that all the $\delta_i$'s are integers.
 Thus, the conclusion of Lemma~\ref{finite tuples} holds
when $\calD(M,\omega)$ is taken to be the set of positive integers
that are smaller than \eqref{bound on delta i}. 
When $[\omega]$ is not rational, Lemma \ref{finite tuples}
follows from the properness property, Lemma \ref{sproper}.
\begin{proof}[Proof of Lemma \ref{finite tuples}]
By Corollary \ref{enough} and Remark \ref{s2triv}, it is enough to show
the lemma in the case that $(M,\omega)$ can be obtained from $\CP^2$ by
a sequence of $\dim H_2(M;\R)-1$ symplectic blow-ups.

Then, by Lemma \ref{sproper}, for any $K > 0$, the set
$\{E \in  H_2(M;\R) \ | \ E \cdot E=-1 \text{ and } 
      0 \leq \omega(E) \leq K\}$ is compact. 
Since $H_2(M;\Z)$ is discrete, its intersection with this compact set,
$${\mathcal{D}}_K = \{E \in H_2(M;\Z) \ | \ 
     E \cdot E=-1 \text{ and }0 \leq \omega(E) \leq K\},$$
is finite.

Let $\Delta$ be a polygon that satisfies the assumptions of 
Lemma \ref{finite tuples}.  Let $\Delta'=\Delta_0,\ldots,\Delta_s=\Delta$ 
be a sequence of Delzant polygons such that each $\Delta_i$ is obtained
from $\Delta_{i-1}$ by a corner chopping of size~$\delta_i$.  
Then $\delta_i$ is the rational length of the new edge $e_i$ of $\Delta_i$
created in the corner chopping.  By part~(1) of Lemma \ref{intersection}, 
the combinatorial self intersection of $e_i$ is $-1$.  
By \S \ref{Z of edges} and Lemma \ref{homomorphism}, 
the composed homomorphism 
$$ \Z[\text{edges of }\Delta_{i}] \rightarrow 
   \Z[\text{edges of }\Delta_{s}]  \rightarrow H_2(M;\Z)$$ 
sends $e_i$ to a class $S \in H_2(M;\Z)$ 
with $\frac{1}{2 \pi} \omega(S)= \delta_i$ and $S \cdot S = -1$.

By \S \ref{eqbound},
$$ \delta_i < \frac{1}{2\pi} K , $$ 
where 
$$ K = 2^{\dim H_2(M;\R)} \int_M \omega \wedge c_1(TM). $$
Also, $0 < \delta_i$ for all $i$.
Thus, the conclusion of Lemma~\ref{finite tuples} holds
when $\calD(M,\omega)$ is taken to be the set 
$\{ \frac{1}{2\pi}\omega(E) \ | \ E \in {\mathcal{D}}_K \}$.
\end{proof}

We give an (elementary) algebraic proof of Lemma \ref{sproper}. 

\begin{proof}[Proof of Lemma \ref{sproper}]
Let $(M,\omega)$ be a symplectic four-manifold that is obtained from
$\tCP$ by $k$ symplectic blow-ups.  The space $H_2(M;\R)$ with the
intersection pairing $\ast \cdot \ast$ is identified  with $R^{1+k}$
with the Minkowski metric, by the basis $L,E_1,\ldots,E_k$ of
$H_2(M)$.  (We denote by $E_1,\ldots,E_k$ the homology classes of the
exceptional divisors obtained by the blow-ups, and by $L$ the homology
class of a line $\CP^{1} \subset M$.)

Let $A \in H_2(M;\R)$ be the Poincar\'e dual to $[\omega]$.  Notice
that $A \cdot A = \omega(A) = \int_{M}{\omega \wedge \omega}$ is
positive, and that $A \cdot C = \omega(C)$ for all $C \in H_2(M;\R)$.
Without loss of generality, $A \cdot A =1$. Complete
$A$ to an orthogonal basis $\B=\{A,B_1,\ldots,B_k\}$ of $\{H_2(M;\R),
\ast \cdot \ast \}$, such that $B_j \cdot B_j = -1$ for all $j$.  The
period map $E \to \omega(E)$ takes $E=x_0 A + \sum_{j=1}^{k} x_j B_j$
to the coefficient $x_0$. The restriction of this map to the
hyperboloid that consists of the classes in $H_2(M;\R)$ with self
intersection $-1$ is proper:  if ${x_0}^2 - \sum_{i=1}^{k} {x_i}^2=-1$
and $x_0$ is bounded, then so are $x_1,\ldots,x_k$.
\end{proof}

\begin{Remark} \labell{sproperg}
Let $(M,\omega)$ be a symplectic four-manifold with $ b_2^+ = 1 $.
The same argument as in the above proof shows that for any $\alpha <0$, 
the restriction of the period map to the set of classes in $H_2(M;\R)$ 
with self intersection $\alpha$ is proper. 
\end{Remark}

\begin{Remark}   \labell{demazure}
For a manifold obtained from $\CP^2$ by eight or fewer blow-ups, the
Poincar\'e dual $C_1$ to the first Chern class $c_1(TM)$ is of positive
self intersection. This implies that the set of homology classes 
$E \in H_2(M;\R)$ that satisfy $E \cdot E = -1$ and $c_1(TM)(E) = 1$ 
is compact. (We can see this by completing $C_1$ to an orthogonal basis
and translating the claim to 
coordinate vectors with respect to this basis.)
Hence the intersection
of this set with the discrete $H_2(M;\Z)$ is finite.  In particular,
there exist only finitely many $E \in H_2(M;\Z)$ of self intersection
$-1$ that can be represented by a symplectically embedded sphere.
See \cite{demazure} or \cite[ch.~IV, \S 26]{manin}.

\end{Remark}

\smallskip

We can now deduce our main theorem.

\begin{proof}[Proof of Theorem \ref{thm:finite}]

Fix a compact connected symplectic four-manifold $(M,\omega)$.

Let $\Delta$ be the moment map image for a toric action on $(M,\omega)$.

By part \eqref{d1} of Proposition \ref{prop:delzant}, $\Delta$ is 
a Delzant polygon.
By parts \eqref{i3} and \eqref{i4} of Lemma \ref{perimeter and area},
\begin{equation} \labell{ho ho}
 \perimeter(\Delta) = \frac{1}{2\pi} \int_M \omega \wedge c_1(TM)
\end{equation}
and
\begin{equation} \labell{ho ho ho}
 \area(\Delta) = \frac{1}{(2\pi)^2} \int_M \frac{1}{2!} \omega \wedge \omega.
\end{equation}

First, suppose that $\dim H_2(M;\R) = 1$. 
By part \eqref{i1} of Lemma \ref{perimeter and area}, 
$\Delta$ has three edges.
By part \eqref{f1} of Lemma \ref{fulton}, $\Delta$ is $\AGL(2,\Z)$-congruent
to a standard Delzant triangle, $\Delta_\lambda$. 
So $\perimeter(\Delta) = 3 \lambda$ and $\area(\Delta) = \half \lambda^2$.
By this and \eqref{ho ho} or \eqref{ho ho ho},
$\lambda$ is uniquely determined by the symplectic manifold $(M,\omega)$.
By part \eqref{d2} of Proposition \ref{prop:delzant},
up to equivalence, there exists only one toric action on $(M,\omega)$.

Now suppose that $\dim H_2(M;\R) \geq 2$.  
By part \eqref{i1} of Lemma \ref{perimeter and area}, 
$\Delta$ has $4+s$ vertices where $s$ is the non-negative integer
\begin{equation} \labell{value of s}
   s = \dim H_2(M;\R) - 2.
\end{equation}
By part \eqref{f2} of Lemma \ref{fulton}, 
after possibly transforming $\Delta$ by an element of $\AGL(2,\Z)$, 
there exist positive numbers $a \geq b > 0$,
an integer $0 \leq k < 2a/b$, and (if $s > 0$) positive numbers
$\delta_1,\ldots,\delta_s > 0$,
such that $\Delta$ can be obtained from the Hirzebruch trapezoid $H_{a,b,k}$
by a sequence of corner choppings of sizes $\delta_1,\ldots,\delta_s$.
By part \eqref{j2} of Lemma \ref{bound},
$$ \perimeter(\Delta) 
   = \perimeter(H_{a,b,k}) - \sum \delta_i
   = 2(a+b) - \sum \delta_i .$$
By part \eqref{j1} of Lemma \ref{bound},
$$ \area(\Delta) 
   = \area(H_{a,b,k}) - \half \sum \delta_i^2
   = ab - \half \sum \delta_i^2 .$$
By this, \eqref{ho ho}, and \eqref{ho ho ho},
the sum and product, $a+b$ and $ab$,
and hence also the numbers $a$ and $b$ themselves,
are determined by the symplectic manifold $(M,\omega)$ 
and the values of $\delta_1,\ldots,\delta_s$.
Because, by Lemma \ref{finite tuples} and by \eqref{value of s}, 
the number of possible tuples 
$(\delta_1,\ldots,\delta_s)$ that can arise in this way is finite,
there are only finitely many possibilities for the values of $a$ and $b$.
Because $k$ is an integer and it satisfies
$$ 0 \leq k < 2a / b ,$$
there are only finitely many possible values for $k$.

We have shown that, after transforming $\Delta$ by an element of $\AGL(2,\Z)$,
there exists a tuple $(a,b,k;\delta_1,\ldots,\delta_s)$
such that $\Delta$ can be obtained from $H_{a,b,k}$ by a sequence
of corner choppings of sizes $\delta_1,\ldots,\delta_s$,
and that the number of tuples $(a,b,k;\delta_1,\ldots,\delta_s)$
that can arise in this way is finite.
For any given polygon, the number of possible corner choppings
of a given size is at most the number of vertices, and hence
it is finite.  It follows that, modulo $\AGL(2,\Z)$-congruence, 
the number of polygons $\Delta$ that can arise as moment map images
of toric actions on $(M,\omega)$ is finite.

By part \eqref{d2} of Proposition \ref{prop:delzant}, the number of 
toric actions on $(M,\omega)$ is finite.
\end{proof}

As explained in the introduction, Theorem \ref{thm:finite} implies that 
for a four dimensional compact connected symplectic manifold $(M,\omega)$ 
with $H^{1}(M;\R)=\{0\}$, the group $\Sympl(M,\omega)$ 
can contain only finitely many non-conjugate two dimensional maximal tori.
We conclude this section with examples of symplectomorphism groups
in which not all maximal tori are two dimensional.
In the first symplectomorphism group of Example \ref{s1t}, 
one maximal torus is two dimensional 
and another maximal torus is one dimensional.  
The second symplectomorphism group contains one dimensional maximal tori 
and does not contain any two dimensional tori.
The third symplectomorphism group does not contain any tori.
The symplectic manifolds in these examples are obtained from $\CP^2$ 
by symplectic blow-ups of equal sizes.  

\begin{Example} \labell{s1t}
Take the $\T^2$-action on $\CP^2$ given by
$(a,b) \cdot [z_0,z_1,z_2] = [z_0,az_1,bz_2]$.
There are exactly three fixed points.
If we restrict to the $S^1$-action given by
$a \cdot [z_0,z_1,z_2] = [z_0, z_1, a z_2]$,
the $S^1$-fixed point set consists of the point 
$[0,0,1]$ and the component $F = \{[z_0,z_1,0]\} \cong \CP^1$.
The moment map for this $S^1$ action attains its minimum along $F$.
See Figure~\ref{fig:CP2} (where $\Phi \colon \CP^2 \to \R^2$
is the $\T^2$ moment map, $i \colon S^1 \to \T^2$ is the inclusion map,
and $i^* \circ \Phi \colon \CP^2 \to \R$ is the $S^1$ moment map).
\begin{figure}[ht]
\setlength{\unitlength}{0.0006in}
\begingroup\makeatletter\ifx\SetFigFont\undefined%
\gdef\SetFigFont#1#2#3#4#5{%
  \reset@font\fontsize{#1}{#2pt}%
  \fontfamily{#3}\fontseries{#4}\fontshape{#5}%
  \selectfont}%
\fi\endgroup%
{\renewcommand{\dashlinestretch}{30}
\begin{picture}(5894,2439)(0,-10)
\put(4650,612){\blacken\ellipse{70}{70}}
\put(3450,612){\blacken\ellipse{70}{70}}
\put(2250,612){\blacken\ellipse{70}{70}}
\put(2250,612){\blacken\ellipse{70}{70}}
\put(1050,612){\blacken\ellipse{70}{70}}
\put(1050,1812){\blacken\ellipse{70}{70}}
\put(1050,612){\blacken\ellipse{70}{70}}
\path(3450,1812)(3450,612)(4650,612)(3450,1812)
\path(525,1212)(600,1137)(525,1062)
\path(2250,1212)(3150,1212)
\path(3075,1287)(3150,1212)(3075,1137)
\path(4650,1212)(5550,1212)
\path(5475,1287)(5550,1212)(5475,1137)
\path(5850,2412)(5850,12)
\put(3450,1812){\blacken\ellipse{70}{70}}
\thicklines
\path(5850,1812)(5850,612)
\path(5850,1812)(5850,612)
\put(4950,1287){\makebox(0,0)[lb]{\smash{{{\SetFigFont{12}{14.4}{\rmdefault}{\mddefault}{\updefault}$\sss i^*$}}}}}
\thinlines
\path(1213,1774)(1227,1767)(1243,1759)
	(1260,1750)(1277,1741)(1295,1729)
	(1315,1717)(1335,1704)(1356,1690)
	(1377,1675)(1400,1659)(1423,1642)
	(1447,1624)(1471,1605)(1496,1586)
	(1520,1566)(1545,1546)(1569,1526)
	(1594,1505)(1618,1484)(1642,1463)
	(1665,1442)(1688,1421)(1710,1400)
	(1732,1379)(1754,1358)(1775,1337)
	(1796,1316)(1817,1294)(1838,1272)
	(1859,1250)(1880,1227)(1901,1204)
	(1922,1180)(1943,1156)(1964,1131)
	(1984,1107)(2004,1082)(2024,1058)
	(2043,1033)(2062,1009)(2080,985)
	(2097,962)(2113,939)(2128,918)
	(2142,897)(2155,877)(2167,857)
	(2179,839)(2188,822)(2197,805)
	(2205,789)(2213,774)(2217,764)
	(2222,753)(2225,743)(2229,733)
	(2231,723)(2234,714)(2235,705)
	(2237,697)(2237,689)(2237,681)
	(2237,674)(2236,667)(2234,660)
	(2232,655)(2229,649)(2226,644)
	(2222,640)(2218,636)(2213,633)
	(2207,630)(2202,628)(2195,626)
	(2188,625)(2181,625)(2173,625)
	(2165,625)(2157,627)(2148,628)
	(2139,631)(2129,633)(2119,637)
	(2109,640)(2098,645)(2088,649)
	(2073,657)(2057,665)(2040,674)
	(2023,683)(2005,695)(1985,707)
	(1965,720)(1944,734)(1923,749)
	(1900,765)(1877,782)(1853,800)
	(1829,819)(1804,838)(1780,858)
	(1755,878)(1731,898)(1706,919)
	(1682,940)(1658,961)(1635,982)
	(1612,1003)(1590,1024)(1568,1045)
	(1546,1066)(1525,1087)(1504,1108)
	(1483,1130)(1462,1152)(1441,1174)
	(1420,1197)(1399,1220)(1378,1244)
	(1357,1268)(1336,1293)(1316,1317)
	(1296,1342)(1276,1366)(1257,1391)
	(1238,1415)(1220,1439)(1203,1462)
	(1187,1485)(1172,1506)(1158,1527)
	(1145,1547)(1133,1567)(1121,1585)
	(1112,1602)(1103,1619)(1095,1635)
	(1088,1649)(1083,1660)(1078,1671)
	(1075,1681)(1071,1691)(1069,1701)
	(1066,1710)(1065,1719)(1063,1727)
	(1063,1735)(1063,1743)(1063,1750)
	(1064,1757)(1066,1764)(1068,1769)
	(1071,1775)(1074,1780)(1078,1784)
	(1082,1788)(1087,1791)(1093,1794)
	(1098,1796)(1105,1798)(1112,1799)
	(1119,1799)(1127,1799)(1135,1799)
	(1143,1797)(1152,1796)(1161,1793)
	(1171,1791)(1181,1787)(1191,1784)
	(1202,1779)(1213,1774)
\path(2150,549)(2137,545)(2123,540)
	(2107,536)(2090,532)(2072,527)
	(2052,523)(2031,519)(2009,515)
	(1986,512)(1961,508)(1936,505)
	(1910,502)(1884,499)(1857,496)
	(1831,494)(1804,492)(1778,491)
	(1752,489)(1726,488)(1700,488)
	(1675,487)(1650,487)(1625,487)
	(1600,488)(1574,488)(1548,489)
	(1522,491)(1496,492)(1469,494)
	(1443,496)(1416,499)(1390,502)
	(1364,505)(1339,508)(1314,512)
	(1291,515)(1269,519)(1248,523)
	(1228,527)(1210,532)(1193,536)
	(1177,540)(1163,545)(1150,549)
	(1142,553)(1134,556)(1126,560)
	(1120,563)(1113,567)(1107,570)
	(1102,574)(1097,578)(1093,582)
	(1089,586)(1085,590)(1083,595)
	(1080,599)(1079,603)(1078,608)
	(1078,612)(1078,616)(1079,621)
	(1080,625)(1083,629)(1085,634)
	(1089,638)(1093,642)(1097,646)
	(1102,650)(1107,654)(1113,657)
	(1120,661)(1126,664)(1134,668)
	(1142,671)(1150,674)(1163,679)
	(1177,684)(1193,688)(1210,692)
	(1228,697)(1248,701)(1269,705)
	(1291,709)(1314,712)(1339,716)
	(1364,719)(1390,722)(1416,725)
	(1443,728)(1469,730)(1496,732)
	(1522,733)(1548,735)(1574,736)
	(1600,736)(1625,737)(1650,737)
	(1675,737)(1700,736)(1726,736)
	(1752,735)(1778,733)(1804,732)
	(1831,730)(1857,728)(1884,725)
	(1910,722)(1936,719)(1961,716)
	(1986,712)(2009,709)(2031,705)
	(2052,701)(2072,697)(2090,692)
	(2107,688)(2123,684)(2137,679)
	(2150,674)(2158,671)(2166,668)
	(2174,664)(2180,661)(2187,657)
	(2193,654)(2198,650)(2203,646)
	(2207,642)(2211,638)(2215,634)
	(2217,629)(2220,625)(2221,621)
	(2222,616)(2222,612)(2222,608)
	(2221,603)(2220,599)(2217,595)
	(2215,590)(2211,586)(2207,582)
	(2203,578)(2198,574)(2193,570)
	(2187,567)(2180,563)(2174,560)
	(2166,556)(2158,553)(2150,549)
\path(1113,1712)(1117,1699)(1122,1685)
	(1126,1669)(1130,1652)(1135,1634)
	(1139,1614)(1143,1593)(1147,1571)
	(1150,1548)(1154,1523)(1157,1498)
	(1160,1472)(1163,1446)(1166,1419)
	(1168,1393)(1170,1366)(1171,1340)
	(1173,1314)(1174,1288)(1174,1262)
	(1175,1237)(1175,1212)(1175,1187)
	(1174,1162)(1174,1136)(1173,1110)
	(1171,1084)(1170,1058)(1168,1031)
	(1166,1005)(1163,978)(1160,952)
	(1157,926)(1154,901)(1150,876)
	(1147,853)(1143,831)(1139,810)
	(1135,790)(1130,772)(1126,755)
	(1122,739)(1117,725)(1113,712)
	(1109,704)(1106,696)(1102,688)
	(1099,682)(1095,675)(1092,669)
	(1088,664)(1084,659)(1080,655)
	(1076,651)(1072,647)(1067,645)
	(1063,642)(1059,641)(1054,640)
	(1050,640)(1046,640)(1041,641)
	(1037,642)(1033,645)(1028,647)
	(1024,651)(1020,655)(1016,659)
	(1012,664)(1008,669)(1005,675)
	(1001,682)(998,688)(994,696)
	(991,704)(988,712)(983,725)
	(978,739)(974,755)(970,772)
	(965,790)(961,810)(957,831)
	(953,853)(950,876)(946,901)
	(943,926)(940,952)(937,978)
	(934,1005)(932,1031)(930,1058)
	(929,1084)(927,1110)(926,1136)
	(926,1162)(925,1187)(925,1212)
	(925,1237)(926,1262)(926,1288)
	(927,1314)(929,1340)(930,1366)
	(932,1393)(934,1419)(937,1446)
	(940,1472)(943,1498)(946,1523)
	(950,1548)(953,1571)(957,1593)
	(961,1614)(965,1634)(970,1652)
	(974,1669)(978,1685)(983,1699)
	(988,1712)(991,1720)(994,1728)
	(998,1736)(1001,1742)(1005,1749)
	(1008,1755)(1012,1760)(1016,1765)
	(1020,1769)(1024,1773)(1028,1777)
	(1033,1779)(1037,1782)(1041,1783)
	(1046,1784)(1050,1784)(1054,1784)
	(1059,1783)(1063,1782)(1067,1779)
	(1072,1777)(1076,1773)(1080,1769)
	(1084,1765)(1088,1760)(1092,1755)
	(1095,1749)(1099,1742)(1102,1736)
	(1106,1728)(1109,1720)(1113,1712)
\path(525,1287)(524,1287)(518,1287)
	(506,1286)(492,1285)(479,1283)
	(468,1281)(459,1278)(450,1274)
	(443,1271)(436,1266)(429,1261)
	(422,1255)(416,1248)(410,1241)
	(406,1233)(403,1226)(401,1219)
	(400,1212)(401,1205)(403,1198)
	(406,1191)(411,1183)(418,1176)
	(425,1169)(434,1163)(443,1158)
	(452,1153)(463,1149)(473,1146)
	(485,1144)(500,1142)(517,1140)
	(537,1139)(558,1138)(578,1138)
	(592,1137)(599,1137)(600,1137)
\path(1530,1092)(1531,1092)(1536,1094)
	(1548,1097)(1566,1102)(1588,1109)
	(1611,1116)(1632,1123)(1650,1131)
	(1666,1138)(1679,1145)(1690,1153)
	(1700,1162)(1709,1172)(1717,1183)
	(1724,1196)(1731,1212)(1739,1230)
	(1746,1251)(1753,1274)(1760,1296)
	(1765,1314)(1768,1326)(1770,1331)(1770,1332)
\dottedline{45}(1770,1332)(1769,1332)(1764,1330)
	(1752,1327)(1734,1322)(1712,1315)
	(1689,1308)(1668,1301)(1650,1293)
	(1634,1286)(1621,1279)(1610,1271)
	(1600,1262)(1591,1252)(1583,1241)
	(1576,1228)(1569,1212)(1561,1194)
	(1554,1173)(1547,1150)(1540,1128)
	(1535,1110)(1532,1098)(1530,1093)(1530,1092)
\path(930,1212)(931,1211)(936,1208)
	(948,1202)(964,1193)(982,1184)
	(999,1176)(1014,1169)(1027,1165)
	(1039,1163)(1050,1162)(1061,1163)
	(1073,1165)(1086,1169)(1101,1176)
	(1118,1184)(1136,1193)(1152,1202)
	(1164,1208)(1169,1211)(1170,1212)
\dottedline{45}(1170,1212)(1169,1213)(1164,1216)
	(1152,1222)(1136,1231)(1118,1240)
	(1101,1248)(1086,1255)(1073,1259)
	(1061,1261)(1050,1262)(1039,1261)
	(1027,1259)(1014,1255)(999,1248)
	(982,1240)(964,1231)(948,1222)
	(936,1216)(931,1213)(930,1212)
\dottedline{45}(1650,492)(1651,493)(1654,498)
	(1660,510)(1669,526)(1678,544)
	(1686,561)(1693,576)(1697,589)
	(1699,601)(1700,612)(1699,623)
	(1697,635)(1693,648)(1686,663)
	(1678,680)(1669,698)(1660,714)
	(1654,726)(1651,731)(1650,732)
\path(1650,732)(1649,731)(1646,726)
	(1640,714)(1631,698)(1622,680)
	(1614,663)(1607,648)(1603,635)
	(1601,623)(1600,612)(1601,601)
	(1603,589)(1607,576)(1614,561)
	(1622,544)(1631,526)(1640,510)
	(1646,498)(1649,493)(1650,492)
\put(0,1137){\makebox(0,0)[lb]{\smash{{{\SetFigFont{12}{14.4}{\rmdefault}{\mddefault}{\updefault}$\sss S^1$}}}}}
\put(450,1812){\makebox(0,0)[lb]{\smash{{{\SetFigFont{12}{14.4}{\rmdefault}{\mddefault}{\updefault}$\sss \CP^2$}}}}}
\put(1800,237){\makebox(0,0)[lb]{\smash{{{\SetFigFont{12}{14.4}{\rmdefault}{\mddefault}{\updefault}$\sss F$}}}}}
\put(2625,1287){\makebox(0,0)[lb]{\smash{{{\SetFigFont{12}{14.4}{\rmdefault}{\mddefault}{\updefault}$\sss \Phi$}}}}}
\end{picture}
}
\caption{}
\labell{fig:CP2}
\end{figure}

Let $T=S^1$ or $T=\T^2$.
For any finite set of $T$-fixed points, we can choose $T$-equivariant
Darboux charts that are centered at these points and are disjoint
from each other; if $\eps > 0$ is sufficiently small,
we can use these Darboux charts to perform $T$-equivariant symplectic
blow-ups of size $\eps$ centered at the chosen fixed points.

\begin{enumerate}
\item
Choose a sufficiently small positive number $\eps$ and perform 
three $S^1$-equivariant symplectic blow-ups of size $\eps$
that are centered at points of $F$.  
Each such a blow-up creates an isolated fixed point whose moment map
value is equal to that of $F$ plus $\eps$.
If a Hamiltonian $S^1$-action on a symplectic four-manifold 
has three or more fixed points with the same moment map value,
the action does not extend to a Hamiltonian $\T^2$-action; 
see \cite{naudin} or \cite[Proposition 5.21]{karshon:periodic}.  
Thus, this $S^1$-action is a maximal torus.
It is a one dimensional maximal torus 
in the symplectomorphism group of a symplectic manifold that is obtained 
from $\CP^2$ by three symplectic blow-ups of size $\eps$.

On the other hand, perform three $\T^2$-equivariant symplectic blow-ups
of size $\eps$, centered at the three $\T^2$-fixed points.
This yields a two dimensional maximal torus in the symplectomorphism group 
of a symplectic manifold that is obtained from $\CP^2$ 
by three symplectic blow-ups of size $\eps$.

By \cite{McDuff-DeformationToIsotopy}, (also see \cite{biran},)
any two symplectic manifolds that are obtained from 
$\CP^2$ by symplectic blow-ups of the same sizes are symplectomorphic.
Therefore, the above constructions give maximal tori in the symplectomorphism
group of the same symplectic manifold; one of these maximal tori is 
one dimensional and the other is two dimensional.

\item

If $\eps$ is sufficiently small, we can perform four 
$S^1$-equivariant symplectic blow-ups centered at points of $F$.
Alternatively, we can perform three $S^1$-equivariant symplectic blow-ups
centered at points of $F$ and one that is centered at the isolated
fixed point $[0,0,1]$.  These two constructions yield two circle actions
that are not equivalent (e.g., the lists of moment map values at their
fixed points are different).  If $1/\eps$ is an integer,
a manifold that is obtained from $\CP^2$ by four or more symplectic blow-ups
of size $\eps$ does not admit a toric action; see \cite[Theorem 4.1 (1)]{kk}.
Thus, the symplectomorphism group of such a symplectic manifold, when $\eps$ is sufficiently 
small and $1/\eps$ is an integer, contains at least two non-conjugate 
one dimensional maximal tori, and it does not contain any two dimensional 
maximal tori.

\item

If $1/\eps$ is an integer and $(k-1)\eps \geq 1$, then a symplectic manifold
that is obtained from $\CP^2$ by $k$ symplectic blow-ups of size $\eps$
does not admit any circle action; see \cite[Theorem 4.1 (2)]{kk}.
It is possible to perform four symplectic blow-ups of $\CP^2$ 
of size $\eps = 1/3$ (see, e.g., \cite{traynor}).
Thus, the resulting symplectic manifold does not admit any circle action. 
\end{enumerate}
\eoe
\end{Example}

\appendix
\section{$J$-holomorphic spheres in symplectic four-manifolds} 
\labell{sec:holomorphic}

The purpose of this appendix is to introduce the non-experts to 
a beautiful work of McDuff that establishes properness for a general
compact symplectic four-manifolds.

We then apply this result, together with a lemma on the existence of certain
symplectically embedded 2-spheres, to give another proof of 
Lemma \ref{finite tuples}.

\begin{Lemma} \labell{finite}
Let $(M,\omega)$ be a closed symplectic four-manifold.
Let ${\mathcal E}$ denote the set of homology classes $E \in H_2(M;\Z)$
that are represented by symplectic exceptional spheres.
Then the map
$E \mapsto \omega(E)$ from ${\mathcal E}$ to $\R$ is proper.
\end{Lemma} 
A \emph{symplectic exceptional sphere} is a symplectically embedded sphere with self intersection $-1$.

Lemma \ref{finite} is proved in \cite[\S 3]{McD:topology}
and \cite[Lemma~3.1]{McDuff-Structure} using Gromov's theory 
of J-holomorphic curves.  We recall relevant parts of this theory. 

Let $(M,\omega)$ be a compact symplectic manifold.
Let  $\J=\J(M,\omega)$ be the space of almost complex structures on $M$ 
that are tamed by $\omega$; see \S\ref{calJ}.  Given $J \in \J$,
a \emph{parametrized $J$-holomorphic sphere} is a map $u \colon \CP^1 \to M$, 
such that $du \colon T\CP^1 \to TM$
satisfies the Cauchy-Riemann equation $du \circ i = J \circ du$.
Such a $u$ represents a homology class in $H_2(M;\Z)$ that we denote $[u]$.
The symplectic area $\omega([u]) := \int_{\CP^1} u^* \omega$ only depends 
on the homology class $[u] \in H_2(M;\Z)$.  
The pull-back $u^* \omega$ is a positively 
oriented area form at each point where $du \neq 0$; it follows that 
the class $[u]$ is non-zero if and only if $u$ is not constant, 
if and only if $\omega([u]) > 0$. A $J$-holomorphic sphere 
is called \emph{simple} if it cannot be factored through a branched 
covering of $\CP^1$.
One similarly defines a holomorphic curve in $(M,J)$ whose domain is a Riemann 
surface other than $\CP^1$.
Gromov, in \cite{gromov}, introduced a notion of ``weak convergence"
of a sequence of holomorphic curves. This notion is preserved under
reparametrization of the curve, and it implies convergence in homology. 
\emph{Gromov's compactness theorem} guarantees that, given a converging
sequence of almost complex structures, a corresponding sequence 
of holomorphic curves with bounded symplectic area has a weakly converging
subsequence.  

The limit under weak convergence might not be a curve; it might
be a ``cusp curve", which is a connected union of holomorphic curves.  
In this paper we don't need the precise definitions 
of cusp curve and of weak convergence.  Readers who are not experienced
with these notions may find this example helpful:

\begin{Example}\labell{ex:cusp curve}
Identify $\CP^1$ with $\C \sqcup \{ \infty \}$.
For each $0 < \lambda < \infty$, consider the holomorphic curve
in $M = \CP^1 \times \CP^1 = 
  (\C \sqcup \{ \infty \}) \times (\C \sqcup \{ \infty \}) $
given by $\{ (w,z) \ | \ w = \lambda z \}$.
This curve can be parametrized by $z \mapsto (\lambda z , z)$
or, alternatively, by $w \mapsto (w,\lambda\Inv w)$.
As $\lambda \to 0$, these curves weakly converge
to the cusp curve 
$ \{ 0 \} \times \CP^1 \cup \CP^1 \times \{ \infty \} $.
\eoe
\end{Example}

\begin{Lemma}\labell{lem:gromov}
Let $\{ J_n \} \subset \J$ be a sequence of almost complex structures
that converges in the $C^\infty$ topology to an almost complex
structure $J_\infty \in \J$.  For each $n$, let $f_n \colon \CP^1 \to M$ 
be a parametrized $J_n$-holomorphic sphere.  Suppose that the set of areas 
$\omega([f_n])$ is bounded.  Then one of the following two possibilities
occurs.
\begin{enumerate}
\item
There exists a $J_\infty$-holomorphic sphere $u \colon \CP^1 \to M$ 
and elements $A_n \in \PSL(2,\C)$ such that a subsequence of the
$f_n \circ A_n $'s converges to $u$ in the $\Cinf$ topology.
In particular, there exist infinitely many $n$'s for which $[f_n] = [u]$.
\item
There exist two or more non-constant simple $J_\infty$-holomorphic
spheres $u_\ell \colon \CP^1 \to M$ and positive integers
$m_\ell$, for $\ell = 1, \ldots, L$, and infinitely many $n$'s for which
$$ [f_n] = \sum_{\ell=1}^L m_\ell [u_\ell] \qquad \text{ in } H_2(M;\Z).$$
\end{enumerate}
\end{Lemma}

\begin{proof}
If there exists a subsequence of the $f_n$'s that, after reparametrization,
converges in the $C^\infty$ topology to a map $u \colon \CP^1 \to M$,
then the map $u$ is $J_\infty$-holomorphic, and the first of the 
two possibilities occurs.

Otherwise, by Gromov's compactness theorem, there exists a subsequence
that weakly converges to a cusp curve.  
See \cite[\S 5]{AL:Sikorav}.
This cusp curve has at least two components; its components are 
non-constant $J_\infty$-holomorphic maps
$v_\ell \colon \CP^1 \to M$;
weak convergence implies convergence in homology, 
so that $[f_n] = \sum_\ell [v_\ell]$ for all but a finite number 
of $f_n$'s in the subsequence.

If $v_\ell$ is not simple then there exists a simple non-constant 
$J_\infty$-holomorphic $u_\ell \colon \CP^1 \to M$ 
and a positive integer $m_\ell$ such that $[v_\ell] = m_\ell [u_\ell]$. 
\end{proof}

\begin{Corollary} \labell{nd finite}
Let $(M,\omega)$ be a closed symplectic manifold.
Let $K$ be a positive real number.  
For each $J \in \J(M,\omega)$, there exist only finitely many 
homology classes $A \in H_2(M;\Z)$ such that
$\omega(A) \leq K$ and such that
the class $A$ is represented by a $J$-holomorphic sphere.
\end{Corollary}

\begin{proof}
Otherwise, there exists a sequence of $J$-holomorphic spheres
$f_n \colon \CP^1 \to M$ such that $\omega([f_n]) \leq K$
for each $n$ and such that the classes $[f_n]$ are all different.
This contradicts Lemma~\ref{lem:gromov}.
\end{proof}


Fix a compact symplectic manifold $(M,\omega)$.
For any class $A\in H_{2}(M;\Z)$, consider the universal moduli space
of simple parametrized holomorphic spheres in the class $A$,
$$ \M(A,\J) = \{ (u,J) \ | \ J \in \J, 
 u \colon \CP^1 \to M \text{ is simple $J$-holomorphic, and } [u]=A \},$$ 
and the projection map
$$ p_{A} \colon \M(A,\J) \to \J. $$
For any positive number $K$, let
$$ \calN_K = \{ A \in H_2(M;\Z) \ | \ 
 A \neq 0 , c_1(TM)(A) \leq 0 , \text{ and } \omega(A) < K \} . $$
The importance of this set lies in the fact that if a homology class $E$
with $\omega(E) \leq K$ and $c_1(TM)(E) \leq 1$ is represented by a
$J$-holomorphic cusp-curve with two or more components
then at least one of these components must lie in a homology class 
in $\calN_K$:

\begin{Lemma}\labell{McD trick}
Let $J \in \J$.  Let
$ u_\ell \colon \CP^1 \to M, \quad \ell = 1, \ldots, L $
be two or more non-constant $J$-holomorphic spheres, 
and let $m_1,\ldots,m_L$ be positive integers.  
Consider the homology class
$ E = \sum m_\ell [u_\ell] $
in $H_2(M;\Z)$.  Suppose that
$ c_1(TM) (E) \leq 1 $.
Let
$ K = \omega (E)  $.
Then there exists $1 \leq \ell \leq L$ such that
$ [u_\ell] \in \calN_K $.
\end{Lemma}

\begin{proof}
By assumption,
\begin{equation} \labell{sum is K}
\sum_{\ell = 1}^L m_\ell \omega([u_\ell]) = K.
\end{equation}
Because the left hand side consists of two or more positive summands,
each summand is smaller than the sum $K$:
\begin{equation} \labell{each lt K}
 \omega([u_\ell]) < K \quad \text{ for all } \ell.
\end{equation}

By assumption,
$$ \sum_{\ell = 1}^L m_\ell c_1(TM)([u_\ell]) \leq 1 .$$
Because this is a positive combination of two or more integers
and the sum is no greater than one, the integers cannot all be
greater than or equal to one, so there must exist at least one $\ell$
such that
\begin{equation} \labell{c1 le 0} 
c_1(TM) ([u_\ell]) \leq 0.
\end{equation}
By this and \eqref{each lt K}, the class $[u_\ell]$ is in $\calN_K$.
\end{proof}

\smallskip

The automorphism group $\PSL(2,\C)$ of $\CP^1$ acts on 
$\M(A,\J)$ by reparametrizations.  
The quotient $\M(A,\J)/\PSL(2,\C)$ is the space of unparametrized
$J$-holomorphic spheres representing $A \in H_2(M;\Z)$.
Let
$$ U_K = \J \ssminus\bigcup_{A \in \calN_K} \image p_A .$$

\begin{Lemma}
Let $E \in H_2(M;\Z)$ be a homology class such that $c_1(TM)(E) \leq 1$.
Let $K = \omega(E)$.  Then the map
$$ p_E\Inv(U_K) / \PSL(2,\C) \to U_K $$
that is induced from the projection map $p_E \colon \M(E,\J) \to \J$
is proper.
\end{Lemma}

\begin{proof}
The lemma follows from Gromov's compactness in the following way.
Let $D \subset U_K$ be a compact subset.
We need to show that $p_E\Inv(D)/\PSL(2,\C)$ is compact.
Because $\M(E,\J)$ is Hausdorff and first countable,
it is enough to show that for every sequence
$\{ (f_n,J_n) \}$ in $p_E\Inv(D)$ there exists a subsequence
that, after reparametrization, has a limit in $p_E\Inv(D)$
in the $C^\infty$ topology.

Take such a sequence, $\{ (f_n,J_n) \} $.
Because $J_n \in D$ and $D$ is compact and contained in $U_K$, 
after passing to a subsequence we may assume that $\{J_n\}$ 
converges to $J_\infty \in U_K$.

Suppose that there exists a subsequence that, after reparametrization,
converges to some $u \colon \CP^1 \to M$ in the $\Cinf$ topology.
Because each $f_n$ is a $J_n$-holomorphic sphere in the class $E$,
the limit $u$ must also be in the class $E$ 
and it must be $J_\infty$-holomorphic.

We need to show that $u$ is simple.  If not,
there is a simple $J_\infty$-holomorphic curve $v$ 
and an integer $m \geq 2$ such that $[u] = m[v]$. 
Because $v$ is simple and $J_\infty$-holomorphic, the pair $(v,J_\infty)$
is in the moduli space $\M(A,\J)$ for $A = [v]$, so 
$J_\infty \in \image p_A$.
The inequalities $c_1(TM)([u]) \leq 1$ and $\omega([u]) = K$
imply that $c_1(TM)([v]) \leq 0$ and $\omega([v]) < K$,
so that $A = [v]$ is in $\calN_K$.  
So $J_\infty \in \image p_A$ cannot be in $U_K$.
Therefore, $u$ is simple.

The pair $(u,J_\infty)$ is then in the moduli space $\M(E,\J)$,
and since $J_\infty \in D$, this pair is in $p_E\Inv(D)$.

Now suppose that there does not exist such a subsequence,
and we would like to arrive at a contradiction.
By Lemma \ref{lem:gromov}, there exist two or more non-constant 
simple $J_\infty$-holomorphic spheres $u_\ell \colon \CP^1 \to M$
and positive integers $m_\ell$ such that $\sum m_\ell [u_\ell] = E$.
Because $\omega(E) = K$ and $c_1(TM)(E) \leq 1$, by Lemma \ref{McD trick} 
at least one of the $u_\ell$'s represents a homology class in $\calN_K$.
This contradicts the fact that the $u_\ell$'s are $J_\infty$-holomorphic
for $J_\infty \in U_K$ and the definition of $U_K$.
\end{proof}

\begin{Corollary} \labell{c1-B}
Let $E$ be a homology class such that $c_1(TM)(E) = 1$
and $\omega(E) \leq K$.  Then $\image p_E \cap U_K$ is closed in $U_K$.
\end{Corollary}

\begin{Lemma}  \labell{c1a}
The subset $U_K \subset \J$ is open.
\end{Lemma}

\begin{proof}
Let $\{ J_n \}$ be a sequence of almost complex structures
that converges in the $C^\infty$ topology to an almost
complex structure $J_\infty$.
We need to show that if each $J_n$ lies in the union
$$ \bigcup_{A \in \calN_K} \image p_A $$
then so does $J_\infty$. For each $n$, let
$f_n \colon \CP^1 \to M$ be a $J_n$-holomorphic sphere such that 
\begin{equation} \labell{inequalities}
 \omega([f_n]) < K \quad \text{ and } \quad
 c_1(TM)([f_n]) \leq 0.
\end{equation}
By Lemma \ref{lem:gromov}, there exist one or more non-constant simple 
$J_\infty$-holomorphic spheres $u_\ell \colon \CP^1 \to M$ 
and positive integers $m_\ell$ such that
\begin{equation} \labell{sum}
\sum_{\ell=1}^L m_\ell [u_\ell] = [f_n]
\end{equation}
for infinitely many $n$'s.  
By Lemma \ref{McD trick} there exists at least one $u_\ell$
whose homology class $[u_\ell]$ is in $\calN_K$.
Because $u_\ell$ is $J_\infty$-holomorphic and $u_\ell$ is simple,
$J_\infty \in \image p_{[u_\ell]}$.
\end{proof}

\smallskip

We have the following consequences of the Sard-Smale theorem 
and the ellipticity of the Cauchy-Riemann equations.

\begin{Lemma} \labell{consequences of sard-smale}
\begin{enumerate}
\item[(a)]
The universal moduli space $\M(A,\J)$ and the space $\J$ of compatible
almost complex structures are Fr\'echet manifolds, and the projection map 
$$p_{A} \colon \M(A,\J) \to \J$$
is a Fredholm map of index 
\begin{equation} \labell{index}
 \dim \ker dp_A - \dim \coker dp_A = 2c_1(TM)(A)+2n, 
\end{equation}
where $2n$ is the dimension of $M$.

\item[(b)]
If $(u,J)$ is a regular value for $p_A$, then for any neighborhood $U$ 
of $(u,J)$ in $\M(A,\J)$, its image, $p_A(U)$, contains a neighborhood 
of $J$ in $\J$.

\item[(c)]
The set of $J$'s that are regular for the map 
$p_A \colon \M(A,\J) \to \J$
for all $A \in H_2(M;\Z)$ is dense in~$\J$.

\item[(d)]
Let $\calA$ be a subset of $H_2(M;\Z)$.  Let $\{ J_t \}_{t \in [0,1]}$
be a $C^1$ simple path in $\J$ whose endpoints are regular values
for $p_A$ for all $A \in \calA$.  Then there exists a $C^1$ perturbation
$\{ \tilde{J}_t \}$ of $J_t$ with the same endpoints
which is transversal to $p_A$ for all $A \in \calA$.

\end{enumerate}
\end{Lemma}

\begin{proof}[Sketch of proof]
For part (a), see \cite[\S 3.1]{MS:JCurves}.

Let $\J^\ell$ denote the space of almost complex structures
of type $C^\ell$.  Let $\M(A,\J^\ell)$ denote the space of pairs
$(u,J)$ where $J \in \J^\ell$ and where $u \colon \CP^1 \to M$
is a simple $J$-holomorphic sphere of type $C^\ell$.
Let $p_A^\ell \colon \M(A,\J^\ell) \to \J^\ell$ denote the projection map.
If $\ell$ is sufficiently large, $\J^\ell$ and $\M(A,\J^\ell)$
are Banach manifolds and $p_A^\ell$ is differentiable;
see \cite[Proposition 3.2.1]{MS:JCurves}.
Moreover, $p_A^\ell$ is a Fredholm map of the same index, $c_1(TM)(A) + 2n$.

The space $\J$ is a dense subset of $\J^\ell$; 
the space $\M(A,\J)$ is a dense subset of $\M(A,\J^\ell)$;
the projection map $p_A$ is the restriction of $p_A^\ell$ to $\M(A,\J)$.
By \cite[Remark 3.2.8]{MS:JCurves}, our notion of ``regular" coincides
with theirs.
A point $(u,J) \in \M(A,\J)$ is regular for $p_A$ if and only if
it is regular for $p_A^\ell$. 

By elliptic regularity, if $J \in \J^\ell$ and $u \colon \CP^1 \to M$
is $J$-holomorphic then $u$ is of type $C^\ell$. 
In particular, if $J \in \J$ is smooth and $u \colon \CP^1 \to M$
is $J$-holomorphic then $u$ is smooth.  See \cite[Prop.~3.1.9]{MS:JCurves}.

Part (b) follows from these facts and the inverse mapping theorem
for Fredholm maps between Banach manifolds.  Explicitly:  because
the differential of $p_A$ at $(u,J)$ is surjective and its kernel
is finite dimensional, the differential of $p_A^\ell$ at $(u,J)$
is also surjective and its kernel is finite dimensional; see 
\cite[last paragraph of \S 3.1]{MS:JCurves}.
In particular, the kernel splits, that is, it has a closed complementary
subspace in $T_{(u,J)} \J(A,\J^\ell)$; see \cite[p.~4]{lang}.
By the inverse function theorem for Banach spaces, there exists
a neighborhood $U^\ell$ of $(u,J)$ in $\M(A,\J^\ell)$ on which $p_A^\ell$
is a projection map; in particular, its image, $p_A(U^\ell)$,
contains an open neighborhood $\Omega^\ell$ of $J$ in $\J^\ell$;
see \cite[chap.~I, \S 5, cor.~2s]{lang}.
After possibly shrinking $U^\ell$, we may assume that its
intersection with $\M(A,\J)$ is contained in $U$.
Since the $C^\ell$ topology on $\J \subset \J^\ell$ is coarser
than the $C^\infty$ topology and $\Omega^\ell$ is open in $\J^\ell$,
the intersection $\Omega^\ell \cap \J$ is an open neighborhood 
of $J$ in $\J$.
Let $J' \in \Omega^\ell \cap \J$.  
Because $J' \in \Omega^\ell$, there exists a $J'$-holomorphic sphere 
$u \colon \CP^1 \to M$ of type $C^\ell$ such that 
$(u,J') \in U^\ell \subset \M(A,\J^\ell)$.
Because $J' \in \J$, by elliptic regularity, $u$ is smooth; 
see \cite[Prop.~3.1.9]{MS:JCurves}.
So $(u,J') \in U^\ell \cap \M(A,\J) \subset U$.
Thus the neighborhood $\Omega^\ell \cap \J$ of $J$ is contained
the image $p_A(U)$.

Part (c) follows from the above facts about the moduli spaces,
together with the Sard-Smale theorem \cite[Theorem (1.3)]{smale} 
for Fredholm maps between Banach manifolds,
see \cite[Theorem 3.1.5 (ii)]{MS:JCurves}, and the fact that 
$H_2(M;\Z)$ is countable.

By \cite[sentence following Definition 3.16]{MS:JCurves},
our notion of ``transversal to $p_A$" coincides with their notion
of a regular path.
Part (d) follows from the above facts about the moduli spaces together 
with the Sard-Smale transversality theorem \cite[Theorem (3.1)]{smale} 
for Fredholm maps between Banach manifolds.
See \cite[Theorem 3.1.7 (ii)]{MS:JCurves}.
\end{proof}

\begin{Lemma}[$C_1$ lemma] \labell{c1}
Let $(M,\omega)$ be a compact symplectic four-manifold. 
The subset $U_K \subset \J$ is dense and path connected.
\end{Lemma}

\begin{proof}
Let $A \in \calN_K$.
Since $A \neq 0$, the $\PSL(2,\C)$-orbits in $\M(A,J)$
are six dimensional submanifolds of the level sets $p\Inv_A(J)$, 
and so $\dim \ker \left. d p_A \right|_{(u,J)} \geq 6$
for all $(u,J) \in \M(A,\J)$.
Since $c_1(TM)(A) \leq 0$ then, by \eqref{index}, 
$\dim \coker dp_A \geq 2$ for all $(u,J) \in \M(A,\J)$.

In particular, $dp_A$ is never onto.
So a point in $\J$ is a regular value for $p_A$ if and only if
it is not in the image of $p_A$.  
The fact that $U_K$ is dense then follows from part (c)
of Lemma \ref{consequences of sard-smale}.

Let $J_0$ and $J_1$ be two points in $U_K$.  Let $\{ J_t \}_{t \in [0,1]}$
be any simple path in $\J$ connecting $J_0$ to $J_1$.  By part (d)
of Lemma \ref{consequences of sard-smale}, there exists a path
$\{ \tilde{J}_t \}_{t \in [0,1]}$ connecting $J_0$ to $J_1$ 
that is transversal to $p_A$ for all $A \in \calN_K$.
Because $\dim \coker dp_A \geq 2$ for all $A \in \calN_K$, 
transversality means that the path is disjoint from $\image p_A$.
So the path lies in $U_K$.
\end{proof}

\begin{Lemma} \labell{adjunction}
Let $(M,\omega)$ be a compact symplectic four-manifold. 
Let $A \in H_2(M;\Z)$ be a homology class that is represented
by an embedded symplectic sphere $C$.  Then for any $J \in \J$ 
and any simple $J$-holomorphic $u \colon \CP^1 \to M$ in the class $A$, 
the map $u$ is an embedding.
\end{Lemma}

\begin{proof}
Let
$ u_0 \colon \CP^1 \to M $
be a symplectic embedding whose image is $C$.
Construct an almost complex structure $J_0 \in \J$ for which 
$u_0$ is a $J_0$-holomorphic sphere.
(Define $J_0|_{TC}$ such that $u_0$ is holomorphic;
extend it to a compatible fiberwise complex structure
on the symplectic vector bundle $TM|_C$; then extend it
to a compatible almost complex structure on $M$.
See \cite[Section 2.6]{MS:intro}.)

By the adjunction inequality, for any $(u,J) \in \M(A,\J)$,
$$ A \cdot A - c_1(TM)(A) + 2 \geq 0, $$
with equality if and only if $u$ is an embedding.
See \cite[Thm.~2.2.1]{AL:McDuff} or \cite[Cor.~E.1.7]{MS:JCurves}.
Applying this to $(u_0,J_0)$, we get that the homology class $A$
satisfies $A \cdot A - c_1(TM)(A) + 2 = 0$.
Applying the adjunction inequality to any other $(u,J) \in \M(A,\J)$,
we get that $u$ is an embedding.
\end{proof}

\begin{Lemma} \labell{repdense}
Let $(M,\omega)$ be a compact symplectic four-manifold.
Let $E \in H_2(M;\Z)$ be a homology class that can be represented
by an embedded symplectic sphere and such that $c_1(TM)(E) = 1$.  Then
\begin{enumerate}
\item
The projection map $ p_E \colon \M(E,\J) \to \J $ is open.
\item
Let $K \geq \omega(E)$. Then, for any $J \in U_K$, 
the class $E$ is represented by an embedded $J$-holomorphic sphere. 
\end{enumerate}
\end{Lemma} 

\begin{proof}
Let $(u,J) \in \M(E,\J)$.  By Lemma \ref{adjunction}, 
the map $u \colon \CP^1 \to M$ is an embedding.
By the Hofer-Lizan-Sikorav regularity criterion \cite{HLS},
if $u$ is an immersed $J$-holomorphic sphere and $c_1([u]) \geq 1$,
then $(u,J)$ is a regular point for the map $p_E \colon \M(E,\J) \to \J$.
By part (b) of Lemma \ref{consequences of sard-smale}
it follows that the map $p_E$ is open. 

Therefore, $\image p_E$ is an open subset of $\J$.
Because, by Lemmas \ref{c1a} and \ref{c1}, $U_K$ is open and dense in $\J$, 
the set $\image p_E \cap U_K$ is a non-empty open subset of $U_K$. 
Because, by Lemma \ref{c1-B}, $\image p_{E} \cap U_K$ is also closed 
in $U_K$, and because, by Lemma \ref{c1}, the set $U_K$ is connected,
we conclude that $\image p_{E} \cap U_K = U_K$. 
So for every $J \in U_K$ there exists a simple $J$-holomorphic sphere
$u \colon \CP^1 \to M$ in the class $E$.
By Lemma \ref{adjunction}, $u$ is an embedding.
\end{proof}

We can now complete the proof of the main result of the appendix.

\begin{proof}[Proof of Lemma \ref{finite}]
Let $K$ be a positive real number.
Suppose that $E_1,E_2,\ldots$ are infinitely many distinct classes in $H_2(M;\Z)$
that satisfy the following conditions.
\begin{enumerate}
\item
The class $E$ can be represented by an embedded symplectic sphere in $M$;
\item
The self intersection $E \cdot E$ is $-1$;
\item
The symplectic area $\omega(E)$ is not larger than $K$.
\end{enumerate}
For each $i$, because $E_i \cdot E_i = -1$ and $E_i$ is represented 
by an embedded symplectic sphere, we have $c_1(TM)(E_i) = 1$.
Pick any $J \in U_K$.  By part (2) of Lemma \ref{repdense}, 
for each~$i$, there exists an embedded $J$-holomorphic sphere
in the class $E_i$.
By Corollary \ref{nd finite}, the $E_i$'s cannot all be distinct.
\end{proof}

To deduce Lemma \ref{finite tuples}, we will also need
the following result,
about the existence of certain symplectic exceptional spheres. 

\begin{Lemma} \labell{repa}
Let $(M,\omega)$ be a compact symplectic four-manifold. Suppose that 
$(M,\omega)$ is obtained by a sequence of 
symplectic blow-ups from a symplectic manifold $(M_0,\omega_0)$.
Let $C_0$ be a symplectic exceptional sphere in $(M_0,\omega_0)$. Then 
$(M,\omega)$ contains a symplectic exceptional sphere $C$ such that $$\int_{C} \omega = \int_{C_0} \omega_0.$$
\end{Lemma}

Because every Delzant polygon $\Delta_0$ can be obtained as the moment
map image of a toric action on a symplectic four-manifold
$(M_0,\omega_0)$ \cite{delzant}, and by Remark \ref{blowup}, we get the
following corollary.

\begin{Corollary} \labell{repb}
Let $(M,\omega)$ be a compact connected symplectic four-manifold.
Let $\Delta$ be the moment map image of a toric action on $(M,\omega)$.
Suppose that $\Delta$ can be obtained from a Delzant polygon $\Delta_0$
by a sequence of corner choppings.  Let $e_0$ be an edge of $\Delta_0$
of combinatorial self intersection $-1$.  Then there exists a
symplectically embedded 2-sphere $C$ in $M$ such that
$$ \frac{1}{2\pi} \int_C \omega = \length (e_0),$$
and the self intersection of $C$ is $-1$.
\end{Corollary}

In fact, a stronger version of Corollary \ref{repb} is true. 
\begin{Lemma} \labell{rep}
Let $(M,\omega)$ be a compact connected symplectic four-manifold.
Let $\Delta$ be the moment map image of a toric action on $(M,\omega)$.
Suppose that $\Delta$ can be obtained from a Delzant polygon $\Delta_0$
by a sequence of corner choppings.  Let $e_0$ be an edge of $\Delta_0$.
Then there exists a symplectically embedded 2-sphere $C$ in $M$ such that
$$ \frac{1}{2\pi} \int_C \omega = \length (e_0), $$
and the self intersection of $C$ is equal to the
combinatorial self intersection of  $e_0$.
\end{Lemma}
This can be proved by ``soft means''; see Appendix \ref{existence}.
Here we show Lemma \ref{repa} using J-holomorphic techniques similar 
to those described above.

\begin{proof}[Sketch of proof of Lemma \ref{repa}] 

Suppose that $(M,\omega)$ 
is obtained from $(M_0,\omega_0)$ by a sequence of 
symplectic
blow-ups: $M_0$, $M_1$, \ldots, $M_s=M$.
Let $\Sigma_0$ be a symplectic exceptional sphere in $(M_0,\omega_0)$.
We claim that there exists a symplectic exceptional sphere 
in $(M,\omega)$ in (the image in $H_2(M)$ of) the class $[\Sigma_0]$.

By induction, it is enough to address the case $s=1$ of a single blow-up.
Let $\iota \colon B^{4} (\delta) \hookrightarrow M_0$
be a 
symplectic embedding defining the blow-up.
By a parametrized version of the ``$C_1$ lemma", 
and because exceptional $J$-holomorphic spheres persist along 
any deformation of regular pairs $\{(\omega_{t},J_{t})\}$, 
see~\cite[Section 3]{McDuff-Structure},
it is enough to show the claim for a deformation equivalent symplectic form
on $M$.  Since any two 
symplectic blow-ups centered at the same fixed point are deformation equivalent,
this implies that we can choose the size $\delta$ arbitrarily small.  
If $\delta$ is small enough, we can find a Hamiltonian isotopy 
$\phi_t \colon M_0 \to M_0$, $t \in [0,1]$,
such that the image of $\phi_1 \circ \iota$ is disjoint from
$\Sigma_0$. Because the embeddings $\iota$ and $\phi_1 \circ \iota$
are isotopic, the symplectic blow-up along $\phi_1 \circ \iota$
gives a manifold $(M',\om')$ symplectomorphic to $(M,\om)$, 
see~\cite[Section 3]{McDuff-Structure}.
Because the symplectic blow-up operation takes place 
in an arbitrarily small neighborhood of the embedded closed ball,
$(M',\om')$, and thus $(M,\omega)$, contains an embedded symplectic sphere
representing the class $[\Sigma_0]$.
\end{proof}

\begin{Remark} \labell{remark existence}
Symplectic exceptional spheres in four-manifolds can be studied
using Taubes-Seiberg-Witten theory of Gromov  
invariants, see~\cite{Li-Liu} and \cite{McDuff-DeformationToIsotopy}.
We can apply this theory to reduce the claim in Lemma \ref{repa} on the
existence
of a symplectic exceptional sphere
to a claim on the existence of an embedded smooth sphere
in the corresponding homology class, and
then deduce the latter claim from the construction of smooth blow-ups. 
\end{Remark}

\begin{noTitle}
We now have holomorphic proof of Lemma \ref{finite tuples}.
Let $(M,\omega)$ be a compact connected four-manifold.
Let $\Delta$ be a moment map image for a toric action on $(M,\omega)$, 
and suppose that $\Delta$ can be obtained from a Delzant polygon $\Delta'$
by a sequence of corner choppings of sizes $\delta_1,\ldots,\delta_s$.  
By \S \ref{eqbound},
$$ \delta_i < \frac{1}{2\pi} K $$ 
for all $i$, where 
$$ K = 2^{\dim H_2(M;\R)} \int_M \omega \wedge c_1(TM). $$ 
By Corollary \ref{repb}, the $\delta_i$'s belong to the set
$$ \calD(M,\omega) := \{ \frac{1}{2\pi} \int_C \omega \ | \ C \subset M
\text{ is a symplectic exceptional sphere } 
\text{ and } \int_C \omega \leq K \}.$$
This set is finite, by Lemma \ref{finite}.
\end{noTitle}

\section{Finding embedded symplectic spheres} 
\labell{existence}

The purpose of this appendix is to prove Lemma \ref{rep}.
This appendix is not included in the published version of this paper. 

We denote the length of an interval $J \subset \R$ by $|J|$.
In what follows we also allow an interval to be degenerate,
that is, to consist of a single point.

Let $\Phi \colon M \to \R^n$ be the moment map for a $\T^n$-action
on $(M,\omega)$.  Let $i \colon S^1 \hookrightarrow \T^n$
be a sub-circle.  The induced map on the dual of the Lie algebras,
$i^* \colon \R^n \to \R$, is a projection in a rational direction;
the moment map for the $i(S^1)$-action is $i^* \circ \Phi$.


\begin{Lemma} \labell{speed1}
Let $(M,\omega)$ be a compact connected $2n$ dimensional symplectic
manifold.  Let $\T^n$ act on $(M,\omega)$ with a moment map $\Phi$
with image $\Delta$.  Let $i \colon S^1 \to \T^n$ be a sub-circle 
with corresponding projection $ i^* \colon \R^n \to \R $.
Let $d$ be an edge of $\Delta$.
Suppose that $i(S^1)$ acts on the two-sphere $\Phi\Inv(d)$
by rotations of speed $k$.  Then
$$ |i^* d| = k \length (d) .$$
\end{Lemma}

\begin{proof}
This is an easy consequence of Lemma \ref{local}.
\end{proof}

\begin{noTitle} \labell{consecutive intervals}
Let $J_1,\ldots,J_n \subset \R$ be (possibly degenerate) closed intervals.
We say that these are \emph{consecutive intervals} if there exist
real numbers $x_0,x_1,\ldots,x_n$ such that, for each $1 \leq i \leq n$,
$J_i$ is the interval with endpoints $x_{i-1}$ and $x_i$.
We say that these intervals are 
\emph{increasing} if $x_0 \leq x_1 \leq \ldots \leq x_n$
and \emph{decreasing} if $x_0 \geq x_1 \geq \ldots \geq x_n$.
\end{noTitle}

\begin{Lemma} \labell{speed2}
Let $i \colon S^1 \to \T^2$ be a sub-circle with corresponding 
projection $i^* \colon \R^2 \to \R$.
Let $\Delta$ be a Delzant polygon.
Let $d$ and $d'$ be two consecutive edges of $\Delta$ 
such that the consecutive intervals $i^*(d)$ and $i^*(d')$ are increasing.
Let $k$ and $k'$ be integers such that $|i^*(d)| = k \length (d)$ 
and $i^*(d') = k' \length(d')$.
Let $\tDelta$ be obtained from $\Delta$ by a corner chopping
at the vertex between $d$ and $d'$, and let $e$ be the resulting new edge.
Then $|i^*(e)|=(k+k')\length(e)$.
\end{Lemma}

\begin{proof}
After $\AGL(2,\Z)$-congruence we may assume that $d$ and $d'$
are contained in the positive coordinate axes 
and $i^*(x,y) = -kx +k'y$. Then
$e = \{ (x,y) \ | \ x \geq 0 , y \geq 0, x+y = \delta \}$,
$\length(e) = \delta$, and $i^*e = [-k\delta,k'\delta]$.
\end{proof}

\begin{Remark}
Let $(M,\omega,\Phi)$ be a toric manifold with moment map image $\Delta$.
The assumption of Lemma \ref{speed2} means that the circle $i(S^1)$
rotates the 2-spheres $\Phi\Inv(d)$ and $\Phi\Inv(d')$ 
with speeds $k$ and $k'$, and that these 2-spheres meet at a fixed point $p$
that is not an isolated minimum or maximum for the moment map
$i^* \circ \Phi$.
By Remark \ref{blowup}, the conclusion of Lemma \ref{speed2} means that
an equivariant symplectic blow-up centered at $p$ creates a 2-sphere
on which the circle $i(S^1)$ acts with speed $k+k'$.
This can also be seen directly: by the local normal form, 
a neighborhood of $p$ is equivariantly symplectomorphic to a neighborhood
of the origin in $\C^2$ with the circle action 
$a \cdot (z,w) = (a^{-k} z, a^{k'}w)$.
Equivariant blow-up creates a $\CP^1$ with homogeneous coordinates $z,w$
on which the circle action 
is $a \cdot [z,w] = [a^{-k} z , a^{k'} w] = [z, a^{k+k'}w]$.
\end{Remark}

Let $\Delta=\Delta_n$ be the moment map image of a toric action 
on $(M,\omega)$; let $\Phi$ denote the moment map. 
Let $\Delta_0,\ldots,\Delta_n$ be a sequence of Delzant polygons
such that each $\Delta_\alpha$ is obtained from $\Delta_{\alpha-1}$
by a corner chopping.  Let $e_0$ be an edge of $\Delta_0$.

Recall the homomorphism~\eqref{injection},
$\Z[\text{edges of }\Delta_{\alpha-1}] 
  \to \Z[\text{edges of }\Delta_{\alpha}]$, discussed in \S\ref{Z of edges}. 
It is easy to show, by induction on $n$, that the image of $e_0$ 
under the composed homomorphism
$$ \Z[\text{edges of }\Delta_{0}] \to \Z[\text{edges of }\Delta_{n}] $$
has the form 
$$ \sum_{\alpha=-j'}^j m_\alpha d_{\alpha} $$
where
$$ d_{-j'}, \ldots, d_{-1}, d_0, d_1, \ldots, d_j $$
are consecutive edges of $\Delta_n$, in counter-clockwise order,
where the coefficients $m_\alpha$ are positive integers,
and where $m_{-j'}$, $m_0$, and $m_j$ are equal to one.

Moreover, $d_0$ is the edge of $\Delta_n$ with the same outward normal
as the edge $e_0$ of $\Delta_0$.
The edge $d_{-j'-1}$ of $\Delta_n$ that immediately precedes $d_{-j'}$ 
in counter-clockwise order has the same outward normal
as the edge $d_-$ of $\Delta_0$ that immediately precedes $e_0$.  Similarly, 
the edge $d_{j+1}$ of $\Delta_n$ that immediately follows $d_{j}$ 
in counter-clockwise order has the same outward normal as the edge $d_+$
of $\Delta_0$ that immediately follows $e_0$.

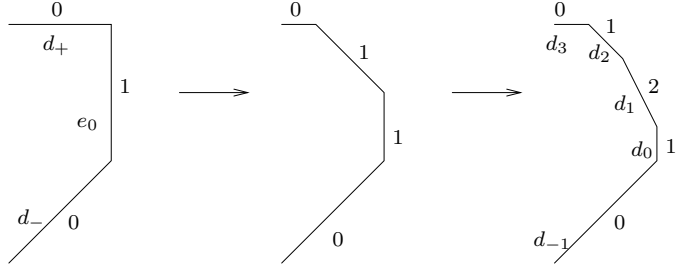
\begin{figure}[ht]
\setlength{\unitlength}{0.0006in}
\begingroup\makeatletter\ifx\SetFigFont\undefined%
\gdef\SetFigFont#1#2#3#4#5{%
  \reset@font\fontsize{#1}{#2pt}%
  \fontfamily{#3}\fontseries{#4}\fontshape{#5}%
  \selectfont}%
\fi\endgroup%
{\renewcommand{\dashlinestretch}{30}
\begin{picture}(5879,2337)(0,-10)
\path(1512,1512)(2112,1512)
\path(1992.000,1482.000)(2112.000,1512.000)(1992.000,1542.000)
\path(2412,2112)(2712,2112)(3312,1512)
	(3312,912)(2412,12)
\path(3912,1512)(4512,1512)
\path(4392.000,1482.000)(4512.000,1512.000)(4392.000,1542.000)
\path(4812,2112)(5112,2112)(5412,1812)
	(5712,1212)(5712,912)(4812,12)
\put(312,1887){\makebox(0,0)[lb]{\smash{{{\SetFigFont{12}{14.4}{\rmdefault}{\mddefault}{\updefault}$\sss d_+$}}}}}
\put(612,1212){\makebox(0,0)[lb]{\smash{{{\SetFigFont{12}{14.4}{\rmdefault}{\mddefault}{\updefault}$\sss e_0$}}}}}
\put(87,362){\makebox(0,0)[lb]{\smash{{{\SetFigFont{12}{14.4}{\rmdefault}{\mddefault}{\updefault}$\sss d_-$}}}}}
\put(387,2187){\makebox(0,0)[lb]{\smash{{{\SetFigFont{12}{14.4}{\rmdefault}{\mddefault}{\updefault}$\sss 0$}}}}}
\put(987,1512){\makebox(0,0)[lb]{\smash{{{\SetFigFont{12}{14.4}{\rmdefault}{\mddefault}{\updefault}$\sss 1$}}}}}
\put(537,312){\makebox(0,0)[lb]{\smash{{{\SetFigFont{12}{14.4}{\rmdefault}{\mddefault}{\updefault}$\sss 0$}}}}}
\put(2487,2187){\makebox(0,0)[lb]{\smash{{{\SetFigFont{12}{14.4}{\rmdefault}{\mddefault}{\updefault}$\sss 0$}}}}}
\path(12,2112)(912,2112)(912,912)(12,12)
\put(3087,1812){\makebox(0,0)[lb]{\smash{{{\SetFigFont{12}{14.4}{\rmdefault}{\mddefault}{\updefault}$\sss 1$}}}}}
\put(4812,2187){\makebox(0,0)[lb]{\smash{{{\SetFigFont{12}{14.4}{\rmdefault}{\mddefault}{\updefault}$\sss 0$}}}}}
\put(3387,1062){\makebox(0,0)[lb]{\smash{{{\SetFigFont{12}{14.4}{\rmdefault}{\mddefault}{\updefault}$\sss 1$}}}}}
\put(2862,162){\makebox(0,0)[lb]{\smash{{{\SetFigFont{12}{14.4}{\rmdefault}{\mddefault}{\updefault}$\sss 0$}}}}}
\put(4627,162){\makebox(0,0)[lb]{\smash{{{\SetFigFont{12}{14.4}{\rmdefault}{\mddefault}{\updefault}$\sss d_{-1}$}}}}}
\put(5487,962){\makebox(0,0)[lb]{\smash{{{\SetFigFont{12}{14.4}{\rmdefault}{\mddefault}{\updefault}$\sss d_0$}}}}}
\put(5337,1362){\makebox(0,0)[lb]{\smash{{{\SetFigFont{12}{14.4}{\rmdefault}{\mddefault}{\updefault}$\sss d_1$}}}}}
\put(5112,1812){\makebox(0,0)[lb]{\smash{{{\SetFigFont{12}{14.4}{\rmdefault}{\mddefault}{\updefault}$\sss d_2$}}}}}
\put(4737,1887){\makebox(0,0)[lb]{\smash{{{\SetFigFont{12}{14.4}{\rmdefault}{\mddefault}{\updefault}$\sss d_3$}}}}}
\put(5337,312){\makebox(0,0)[lb]{\smash{{{\SetFigFont{12}{14.4}{\rmdefault}{\mddefault}{\updefault}$\sss 0$}}}}}
\put(5787,987){\makebox(0,0)[lb]{\smash{{{\SetFigFont{12}{14.4}{\rmdefault}{\mddefault}{\updefault}$\sss 1$}}}}}
\put(5637,1512){\makebox(0,0)[lb]{\smash{{{\SetFigFont{12}{14.4}{\rmdefault}{\mddefault}{\updefault}$\sss 2$}}}}}
\put(5262,2037){\makebox(0,0)[lb]{\smash{{{\SetFigFont{12}{14.4}{\rmdefault}{\mddefault}{\updefault}$\sss 1$}}}}}
\end{picture}
}

%
%
%
%
%
\caption{Coefficients of the image of $e_0$ in $\Z[\text{edges of }\Delta_i]$
for $i=0,1,2$}
\labell{fig:induction}
\end{figure}

Now consider the consecutive edges $d_-,e_0,d_+$ of $\Delta_0$.
Let $i \colon S^1 \to \T^2$ be the unique sub-circle
such that $i^* d_+$ is a single point and oriented so that 
the three consecutive intervals $i^* d_-$, $i^* e_0, i^* d_+$
are increasing.  Then $|i^* e_0| = \length (e_0)$ and $|i^* d_+| = 0$.
Similarly, let $i' \colon S^1 \to \T^2$ be the unique sub-circle
such that ${i}'^* d_-$ is a single point and oriented so that 
the three consecutive intervals ${i'}^* d_-$, ${i'}^* e_0, {i'}^* d_+$
are increasing.  Then $|{i'}^* e_0| = \length (e_0)$ and $|{i'}^* d_-| = 0$.
Without loss of generality, we may assume that the origin is contained
in $e_0$, so that
the restriction to $e_0$ of $i$ coincides with that of $i'$.
\begin{figure}[ht]
\setlength{\unitlength}{0.0006in}
\begingroup\makeatletter\ifx\SetFigFont\undefined%
\gdef\SetFigFont#1#2#3#4#5{%
  \reset@font\fontsize{#1}{#2pt}%
  \fontfamily{#3}\fontseries{#4}\fontshape{#5}%
  \selectfont}%
\fi\endgroup%
{\renewcommand{\dashlinestretch}{30}
\begin{picture}(1710,2739)(0,-10)
\put(1062,162){\makebox(0,0)[lb]{\smash{{{\SetFigFont{12}{14.4}{\rmdefault}{\mddefault}{\updefault}${i'}^*$}}}}}
\path(837,2412)(1287,2412)
\path(1167.000,2382.000)(1287.000,2412.000)(1167.000,2442.000)
\path(837,237)(1287,687)
\path(1223.360,580.934)(1287.000,687.000)(1180.934,623.360)
\put(912,2487){\makebox(0,0)[lb]{\smash{{{\SetFigFont{12}{14.4}{\rmdefault}{\mddefault}{\updefault}$\sss i^*$}}}}}
\put(1587,2487){\makebox(0,0)[lb]{\smash{{{\SetFigFont{12}{14.4}{\rmdefault}{\mddefault}{\updefault}$\sss \R$}}}}}
\path(12,2262)(1512,2262)(1512,1062)
	(612,162)(312,162)(12,462)(12,2262)
\dottedline{45}(1512,2712)(1512,12)
\end{picture}
}
\caption{}
\labell{fig:iiprime}
\end{figure}

For each $-(j'+1) \leq \alpha \leq j+1$, consider the two-sphere
$C_\alpha = \Phi\Inv(d_\alpha)$.
The circle $i(S^1)$ rotates the two-spheres $C_0$ and $C_j$ 
with speed 1 and fixes the two-sphere $C_{j+1}$; its moment map,
$i^* \circ \Phi \colon M \to \R$, takes its maximum along the sphere
$C_{j+1}$.  Similarly, The circle $i'(S^1)$ rotates the two-spheres $C_0$ 
and $C_{-j'}$ with speed 1 and fixes the two-sphere $C_{-(j'+1)}$;
its moment map, $(i')^* \circ \Phi$, takes its minimum along the sphere
$C_{-j'-1}$.  On the sphere $C_0$, the circle actions $i(S^1)$ 
and $i'(S^1)$ coincide.

By iterations of Lemma \ref{speed2}, for each $0 \leq \alpha \leq j$,
$$ |i^* d_\alpha| = m_{\alpha} \length (d_\alpha) ,$$
and for each $0 \leq \alpha \leq j'$,
$$ |{i'}^* d_{-\alpha}| = m_{{-\alpha}} \length (d_{-\alpha}). $$

For each $0 \leq \alpha \leq j$ the circle $i(S^1)$ rotates the sphere 
$C_\alpha$ with speed $m_\alpha$, and for each $0 \leq \alpha \leq j'$
the circle ${i'}(S^1)$ rotates the sphere $C_{-\alpha}$ 
with speed $m_{-\alpha}$; this follows from Lemma \ref{speed1}.

Let $\gamma \colon [-1,1] \to M$ be any path for which 
$\gamma(-1) \in C_{-j'-1}$, $\gamma(1) \in C_{j+1}$,
and $\gamma(t) \in C_0$ for $- \eps < t < \eps$.
Consider the map 
$$ \tgamma \colon S^1 \times [-1,1] \to M $$
that is obtained from $\gamma$ by ``sweeping" by the circle actions:
\begin{equation} \labell{sweep}
 \tgamma(a,t) = \begin{cases}
 i(a) \cdot \gamma(t) & \text{ for } t \in [-\eps,1] \\
 i'(a) \cdot \gamma(t) & \text{ for } t \in [-1,\eps] .
\end{cases} 
\end{equation}
This is well defined for $t \in (-\eps,\eps)$ because $\gamma(t) \in C_0$
for such $t$ and because the actions $i(S^1)$ and $i'(S^1)$ coincide
on $C_0$.  
Because $\tgamma$ is constant on $S^1 \times \{ -1 \}$ and on 
$S^1 \times \{ 1 \}$, it represents a homology class, $[\tgamma]$.
Because $M$ is simply connected, this homology class
is independent of the choice of path $\gamma$ with the above properties.

Starting from a path $\gamma$ that traces a chain of meridians in the spheres
$$ C_{-j'}, \ldots, C_0, \ldots, C_j ,$$
the resulting class $[\tgamma]$ is equal to 
$\sum_{\alpha=-j'}^j m_\alpha [C_\alpha]$.
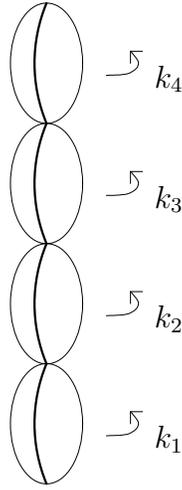
\begin{figure}[ht]
\setlength{\unitlength}{0.00083333in}
\begingroup\makeatletter\ifx\SetFigFont\undefined%
\gdef\SetFigFont#1#2#3#4#5{%
  \reset@font\fontsize{#1}{#2pt}%
  \fontfamily{#3}\fontseries{#4}\fontshape{#5}%
  \selectfont}%
\fi\endgroup%
{\renewcommand{\dashlinestretch}{30}
\begin{picture}(1092,3075)(0,-10)
\path(608,1830)(609,1830)(616,1830)
	(630,1831)(650,1831)(671,1832)
	(691,1833)(708,1835)(723,1837)
	(735,1839)(746,1842)(756,1846)
	(765,1851)(774,1856)(783,1862)
	(790,1869)(797,1876)(802,1884)
	(805,1891)(807,1898)(808,1905)
	(807,1914)(803,1923)(796,1934)
	(786,1947)(775,1961)(765,1973)
	(759,1979)(758,1980)
\path(758,1155)(758,1230)(833,1230)
\path(608,1080)(609,1080)(616,1080)
	(630,1081)(650,1081)(671,1082)
	(691,1083)(708,1085)(723,1087)
	(735,1089)(746,1092)(756,1096)
	(765,1101)(774,1106)(783,1112)
	(790,1119)(797,1126)(802,1134)
	(805,1141)(807,1148)(808,1155)
	(807,1164)(803,1173)(796,1184)
	(786,1197)(775,1211)(765,1223)
	(759,1229)(758,1230)
\path(758,405)(758,480)(833,480)
\path(608,330)(609,330)(616,330)
	(630,331)(650,331)(671,332)
	(691,333)(708,335)(723,337)
	(735,339)(746,342)(756,346)
	(765,351)(774,356)(783,362)
	(790,369)(797,376)(802,384)
	(805,391)(807,398)(808,405)
	(807,414)(803,423)(796,434)
	(786,447)(775,461)(765,473)
	(759,479)(758,480)
\path(758,2655)(758,2730)(833,2730)
\path(608,2580)(609,2580)(616,2580)
	(630,2581)(650,2581)(671,2582)
	(691,2583)(708,2585)(723,2587)
	(735,2589)(746,2592)(756,2596)
	(765,2601)(774,2606)(783,2612)
	(790,2619)(797,2626)(802,2634)
	(805,2641)(807,2648)(808,2655)
	(807,2664)(803,2673)(796,2684)
	(786,2697)(775,2711)(765,2723)
	(759,2729)(758,2730)
\thicklines
\put(1133.000,2655.000){\arc{1950.000}{2.7468}{3.5364}}
\put(1133.000,1905.000){\arc{1950.000}{2.7468}{3.5364}}
\thinlines
\path(758,1905)(758,1980)(833,1980)
\thicklines
\put(1133.000,405.000){\arc{1950.000}{2.7468}{3.5364}}
\put(908,255){\makebox(0,0)[lb]{\smash{{{\SetFigFont{12}{14.4}{\rmdefault}{\mddefault}{\updefault}$k_1$}}}}}
\put(1133.000,1155.000){\arc{1950.000}{2.7468}{3.5364}}
\thinlines
\put(233,2655){\ellipse{450}{750}}
\put(233,1905){\ellipse{450}{750}}
\put(233,1155){\ellipse{450}{750}}
\put(233,405){\ellipse{450}{750}}
\put(908,1005){\makebox(0,0)[lb]{\smash{{{\SetFigFont{12}{14.4}{\rmdefault}{\mddefault}{\updefault}$k_2$}}}}}
\put(908,2505){\makebox(0,0)[lb]{\smash{{{\SetFigFont{12}{14.4}{\rmdefault}{\mddefault}{\updefault}$k_4$}}}}}
\put(908,1755){\makebox(0,0)[lb]{\smash{{{\SetFigFont{12}{14.4}{\rmdefault}{\mddefault}{\updefault}$k_3$}}}}}
\end{picture}
}
\caption{sweeping a chain of meridians}
\label{fig:chain}
\end{figure}
See Figure~\ref{fig:chain}.
By Lemma \ref{homomorphism} and since the homomorphism \eqref{injection} 
respects the length
functional and the combinatorial intersection pairing, 
the pairing of the class $[\tgamma]$ with $\omega$
is equal to $2\pi \length(e_0)$, and the self intersection of the class
$[\tgamma]$ is equal to the self intersection of $e_0$.

To prove Lemma \ref{rep}, it is enough to find a path $\gamma$
with the above properties such that $\tgamma$ defines a symplectically 
embedded sphere.


There exists a $\T^2$-invariant complex structure $J$ on $M$ 
that is compatible with the symplectic form $\omega$ (see \cite{naudin}).
The spheres $C_\alpha$ are invariant under the action of the 
complexified torus $(\C^*)^2$.
If $\gamma(t)$ is a gradient trajectory of $i^* \circ \Phi$
with respect to the metric defined by $J,\omega$
then $\tgamma(a,t)$, defined in \eqref{sweep},
is a $\C^*$-orbit for the complexification of $i(S^1)$.

Let $U$ and $U'$ be $\T^2$-invariant
open subsets of $M$ whose closures are disjoint, such that
$C_1 \cup \ldots \cup C_j \subset U$ and
$C_{-1} \cup \ldots \cup C_{-j'} \subset U'$.
Take an orbit in $C_0 \cap U$.
By \cite[Lemma 3.6]{karshon:periodic}, we can perturb the metric
on a small neighborhood of this orbit to obtain a compatible metric
whose restriction to $U$ is $i(S^1)$-invariant
with the property that an upward gradient trajectory of $i^* \circ \Phi$
starting from a point on $C_0 \cap U'$ approaches, in the limit, 
a point on $C_{j+1}$.
Similarly, we can perturb the metric in $U'$ to obtain a compatible metric 
whose restriction to $U'$ is $i'(S^1)$-invariant such that a downward
gradient trajectory of $i' \circ \Phi$ starting from a point in $C_0 \cap U$
approaches, in the limit, a point on $C_{-j'-1}$.

Take such a gradient trajectory.  After reparametrization, it extends 
to a continuous map
$$ \gamma \colon [-1,1] \to M$$
that is smooth on the open interval $(-1,1)$,
such that $\gamma(t) \in C_0$ for $t \in (-\eps,\eps)$ for some $\eps > 0$, 
such that $\frac{d}{dt} \gamma(t)$ is a positive multiple 
of $\text{grad} (i^* \circ \Phi)$ for $0 \leq t < 1$
and is a positive multiple of $\text{grad} ({i'}^* \circ \Phi)$ 
for $-1 < t \leq 0$,
and such that $\gamma(1) \in C_{j+1}$ and $\gamma(-1) \in C_{-j'-1}$.

For $(a,t)$ with $-1 < t < 1$, $d\tgamma|_{(a,t)} \colon \R^2 \to TM$
is an inclusion onto a symplectic subspace.
A neighborhood of $\gamma(1)$ in $M$ is equivariantly biholomorphic 
to a neighborhood of the origin in $\C^2$ where $i(S^1)$ acts on the 
first coordinate.   In this neighborhood, the closures 
of the $\C^*$-orbits are smooth.
A similar argument holds for $t$ near $-1$. 
Hence, the image of $\tgamma$ 
is a symplectically embedded two-sphere $C$ in $M$.


This completes the proof of Lemma \ref{rep}.

\end{document}